\documentclass[11pt]{amsart}

\usepackage{amsfonts}
\usepackage{amsmath}
\usepackage{amsthm}
\usepackage{mathrsfs}
\usepackage{enumerate}
\usepackage{graphicx}
\usepackage{caption}
\usepackage{subcaption}
\usepackage{caption}
\usepackage[square,sort,comma,numbers]{natbib}
\usepackage{color}
\usepackage{hyperref}
\usepackage{setspace}
\usepackage{natbib}
\usepackage[T1]{fontenc}
\usepackage{times}
\usepackage[left=2cm,right=2cm,top=2.5cm,bottom=2.5cm]{geometry}

\usepackage{bigints}

\newcommand{\EXCLUDE}[1]{}

\newcommand{\remove}[1]{}
\newcommand{\beq}{\begin{eqnarray}}
	\newcommand{\eeq}{\end{eqnarray}}
\newcommand{\beqq}{\begin{eqnarray*}}
	\newcommand{\eeqq}{\end{eqnarray*}}%\pagestyle{plain}
\def\:{:\,}

\newtheorem{thm}{Theorem}[section]

\newtheorem{lem}[thm]{Lemma}

\numberwithin{equation}{section}

\def\rar{\rightarrow}

\newcommand{\be}{\beta}

\newcommand{\ep}{\epsilon}
\newcommand{\al}{\alpha}

\newcommand{\del}{\delta}
\newcommand{\lam}{\lambda}
\newcommand{\Lam}{\Lambda}

\def\th{\theta}

\def\bR{\mathbb{R}}

\def\bZ{\mathbb{Z}}

\def\bN{\mathbb{N}}

\def\bX{\mathbb{X}}

\newcommand{\cP}{\mathcal{P}}

\newcommand{\cY}{{\mathcal Y}}

\def\1{\mathbf{1}}

\def\lab{\label}

\def\nn{\nonumber}
\def\f{\frac}

\begin{document}
	
	%%%%%%%%%%%%%%%%%%%%%%%%%%%%%%%%%%%%FRONTPAGE%%%%%%%%%%%%%%%%%%%%%%%%%%%%%%%%%%%%%
	%\bibliographystyle{abbrv}
	
	% "Title of the paper"
	\title[Connectivity in a Scale-free network]{\large{Poisson Approximation and Connectivity in a Scale-free Random Connection Model}}

\author{Srikanth K. Iyer}
\address{(SKI) Department of Mathematics\\
		Indian Institute of Science\\
	       Bangalore, India.\\
	       Email-\textit{srikiyer@gmail.com}}
%\email{skiyer@iisc.ac.in}

\author{Sanjoy Kr. Jhawar}
\address{(SKJ) Department of Mathematics\\
	Indian Institute of Science \\
	Bangalore, India.\\
	Email-\textit{sanjayjhawar@iisc.ac.in}}
%\email{sanjayjhawar@iisc.ac.in}

\noindent  \thanks{ SKI's research was supported in part from Matrics grant from SERB and DST-CAS.
 SKJ's research was supported by DST-INSPIRE Fellowship. Corresponding author email: srikiyer@gmail.com }
 \maketitle
\begin{abstract}
 We study an inhomogeneous random connection model in the connectivity regime. The vertex set of the graph is a homogeneous Poisson point process $\mathcal{P}_s$ of intensity $s>0$ on the unit cube $S=\left(-\frac{1}{2},\frac{1}{2}\right]^{d},$ $d \geq 2$ . Each vertex is endowed with an independent random weight distributed as $W$, where $P(W>w)=w^{-\beta}1_{[1,\infty)}(w)$, $\beta>0$. Given the vertex set and the weights an edge exists between $x,y\in \mathcal{P}_s$ with probability 
$\left(1 - \exp\left( - \frac{\eta W_xW_y}{\left(d(x,y)/r\right)^{\alpha}} \right)\right),$
independent of everything else, where $\eta, \alpha > 0$, $d(\cdot, \cdot)$ is the toroidal metric on $S$ and $r > 0$ is a scaling parameter.  We derive conditions on $\alpha, \beta$ such that under the scaling
$r_s(\xi)^d= \frac{1}{c_0 s} \left( \log s +(k-1) \log\log s +\xi+\log\left(\frac{\alpha\beta}{k!d} \right)\right),$
$\xi \in \bR$, the number of vertices  of degree $k$ converges in total variation distance to a Poisson random variable with mean $e^{-\xi}$ as $s \to \infty$, where $c_0$ is an explicitly specified constant that depends on $\al, \be, d$ and $\eta$ but not  on $k$. In particular, for $k=0$ we obtain the regime in which the number of isolated nodes stabilizes, a precursor to establishing a threshold for connectivity. We also derive a sufficient condition for the graph to be connected with high probability for large $s$. The Poisson approximation result is derived using the Stein's method. 
\end{abstract}


\section{Introduction}

Social, financial and other networks such as the internet and wireless networks have been objects of much interest among researchers and practitioners in various fields in recent years. This is largely due to following stylized features observed in empirical data (See \citep{Newman2566}, and Section 1.3, in \citep{Durrett2006}). \textit{Small world effect:} Typically the number of `links' required to connect two distant vertices is very small (see \citep{Wattts6Degrees}). \textit{Clustering property:} Linked vertices have mutual connections. \textit{Heavy-tailed degree distribution:} It has been observed that the degree distribution is heavy-tailed with tail parameter between $1$ and $2$ (see \citep{Durrett2006}).

These real-life networks naturally possess long range connections among the vertices. The homogeneous long range percolation model on $\bZ$ was first introduced by Zhang in \citep{Zhang83}. The model can be extended to a random graph with vertex set $\bZ^d$. Vertices $x, y$ are connected by an edge with probability proportional to $\eta |x-y|^{-\al}$ as $|x-y| \to \infty$ for some $\eta, \al>0$. Since nearby points are connected with higher probability, this form of connection probability leads to the graph having the clustering property. For certain range of values of the parameter $\al$ a small world effect has also been observed \citep{Biskup2004}. The continuum version of the above graph called the random connection model was introduced in \citep{Penrose1991}. Such models have found applications in the study of wireless communication networks, models for spread of epidemic and interactions among molecules \citep{Franceschetti2007}. The vertex set of the graph is a homogeneous Poisson point process $\cP_{\lam}$ of intensity $\lam > 0$ in $\bR^d$. $x,y \in \cP_{\lam}$ are connected by an edge with probability $g(x-y)$ where $g:\bR^d \to [0,1]$. It was shown that a non-trivial phase transition occurs if and only if $g$ is integrable. This set-up where the expected degree is finite is referred to as the thermodynamic regime. Of interest are a non-trivial threshold for percolation, the degree distribution and the graph distance.

The random graph models described above do not exhibit a heavy-tailed degree distribution. To overcome this, Deijfen et. al. in \citep{Deijfen2013} proposed an inhomogeneous version of this long-range percolation model on $\bZ^d$ and called it the scale-free percolation model. The inhomogeneity was introduced by assigning independent and identically distributed weights $W_x$  at each vertex $x \in \bZ^d$ representing of the importance of a vertex. Given any two points $x,y \in \bZ^d$ and corresponding weights $W_x, W_y$, the probability that there is an edge between these two points equals $1 - \exp{\left(- \f{\eta W_x W_y}{|x-y|^{\al}}\right)}$ where $\eta, \al > 0$ and $P(W_x > w) = w^{-\beta}1_{[1,\infty)}(w)$ for some $\beta > 0$. For the graph to be non-trivial (finite degrees) one must have $\min\{\al, \al \beta\} > d$. This model is studied in great detail in \citep{Deprez2015}. When $\al \beta < 2d$ the degree distribution is heavy tailed. Surprisingly, under this condition the graph also shows an ultra small world effect, that is, the graph distance between two far away points grows doubly logarithmically in the Euclidean distance. A non-trivial phase transition occurs in the model if $\al \beta > 2d$. A non-trivial phase transition for the above model refers to the existence of a critical value $\eta_c \in (0, \infty)$ such that for any $\eta < \eta_c$ all components in the graph are finite and for all $\eta > \eta_c$, there is, with probability one, an infinite component in the graph. In addition, if $\al \in (d,2d)$ the graph displays a small-world effect. Thus this model is rich enough to exhibit all the stylized features of real-world networks. The continuity of the percolation function which was conjectured in \citep{Deijfen2013} was proved in \citep{Deprez2015}. The above results on $\bZ^d$ were extended to the continuum in \citep{Deprez2018}. 

%The vertices in the graph are distributed according to a homogeneous Poisson point process in $\bR^d$. 
%In the continuum setting the Campbell-Mecke formula enables one to obtain explicit formulas for the graphical quantities of interest. It is this model that we shall study in the connectivity regime. 
%The graph such as the one above where the weights are degenerate belong to a class of geometric random graphs called the random connection models, which we shall describe in brief.

The graph in the thermodynamic regime is however far from being connected. To obtain connectivity, one needs to scale the connection function suitably. As is the case with the Erdos-Renyi random graph, though much harder to prove, the main obstacle to connectivity in certain random geometric graphs is the presence of isolated nodes (see Chapter 13, \citep{Penrose2003}, \citep{Penrose2016}). A scaling in the connectivity regime under which the number of isolated vertices converges to a Poisson distribution for the random connection model is obtained in \citep{MaoAnderson2011}. Suppose that the graph is connected with probability approaching one in the absence of isolated nodes. The Poisson convergence result then yields the asymptotic probability that the graph is connected. Such a result is available only in some restricted cases for the random connection model. The soft random geometric graph was considered by Penrose in \citep{Penrose2016}. In this model each pair of $n$ points distributed uniformly in the unit square are connected with probability $p$ provided the inter-point distance between them is at most $r$. Conditions for a Poisson approximation for the number of isolated nodes and the asymptotic probability of connectivity are derived for arbitrary sequences of parameters $(p_n,r_n)$. The former result is extended to a larger class of connection functions in higher dimensions. A sufficient condition for the random connection model to be connected asymptotically with high probability is derived in \citep{Iyer2018}.

Another model of interest with a different type of inhomogeneity is a general random connection model studied by \citep{Penrose2018} in the connectivity regime. Non-uniformity in this graph comes from two sources: the vertices of the graph form a Poisson point process with intensity measure $s \mu$ where $\mu$ is a probability measure and $s > 0$ is a parameter. Given a realization of the vertices, any two pair of vertices are connected with probability $\phi_s(x,y)$, where $\phi_s : \bX \times \bX \to [0, 1]$ is a symmetric function. Of interest are the number of nodes of a certain degree and the number of components of a given size. It is shown that the number of vertices of a fixed degree converge (as $s \to \infty$) to a Poisson distribution when the expected number of such vertices stabilize. A key assumption under which the result is proved is that the connection function satisfies $\max_{x, y} \phi_s(x,y) < 1 - \ep$ for some $\ep > 0$, that is, it stays bounded away from one. Thus it does not include the simple Gilbert's disk model where an edge exists between any two pair of nodes that are within a specified distance from each other. That the above condition does not hold for our model complicates the computations.

We shall study, for the inhomogeneous random graph considered in \citep{Deprez2018}, the asymptotic behavior of the number of isolated nodes and vertices of arbitrary degree in the connectivity regime. As is often the case, it is convenient to study the connectivity problem on the unit cube. The vertex set is assumed to be distributed according a homogeneous Poisson point process $\cP_s$ of intensity $s$ on the unit cube. To avoid boundary effects, we equip the space with the toroidal metric $d(\cdot, \cdot)$ (to be specified below). The connection function is given by 
\begin{equation}
g\left(x,W_x; y, W_y\right)=1-\exp\left(-\f{\eta W_x W_y}{\left(d(x,y)\right)^{\al}}\right),
\lab{eq:connection_function}
\end{equation}
where $\eta, \al > 0$ and $P(W_x > w) = w^{-\beta}1_{[1,\infty)}(w)$ for some $\beta > 0$. 
To stabilize the expected number of isolated nodes and vertices of fixed degree we must scale the connection function by a parameter depending on the intensity that we denote by $r_s$. The modified connection function will be denoted by
\begin{equation}
	g_s\left(x,W_x; y, W_y\right) = 1-\exp\left(-\f{\eta W_x W_y}{\left(\f{d(x,y)}{r_s}\right)^{\al}}\right).
	\lab{eq:connection_function_s}
\end{equation}
We derive an explicit expression for this scaling function and show that the number of vertices of arbitrary degree converges to a Poisson distribution under certain condition on the parameters. 
The Poisson convergence results are proved using the Stein's method and thus give a bound on the rate of convergence in the total variation distance. We adapt Theorem 3.1 from \citep{Penrose2018} to account for the inhomogeneity arising from the random weights associated with the vertices of the graph. Much of the effort in proving this result lies in overcoming the challenge posed by the fact that the connection function given by (\ref{eq:connection_function_s}) can take values arbitrarily close to one and the presence of weights on the vertices. The connectivity problem for geometric random graphs is, in general hard. Having derived the Poisson convergence for the number of isolated nodes, one would like to show that the graph is connected with high probability whenever the isolated nodes vanish. We derive a sufficient condition for the graph to be connected with high probability under mild conditions on the parameters.

\section{The inhomogeneous random connection model}

We now provide a precise definition of the model that is the object of our study. 
Consider the unit cube $S=\big(-\f{1}{2},\f{1}{2}\big]^d$. We shall ignore ``edge effects'' by equipping $S$ with the toroidal metric $d(\cdot, \cdot) : S \times S \rar \bR^+ \cup \{0\}$ defined by $d(x,y)=  \inf\{ ||x-y-z|| : z \in \bZ^d   \}$, where $|| \cdot ||$ is the Euclidean norm. Let $\cP_s$ be a Poisson point process with intensity $s$ on $S$. Consider the random graph with vertex set $\cP_s$. In this model the edge probabilities depends on the distance as well as the weights at the vertices. To each $x \in \cP_s$ we associate independent random weights with probability distribution satisfying $ P(W > w)=w^{-\be}1_{[1, \infty)}(w)$, $\be>0$. The probability that there is an edge between vertices (located at) $x,y \in \cP_s$ is given by  (\ref{eq:connection_function_s})
independent of everything else for fixed $\al>0, \eta>0$. We denote the resulting random graph by $G(\cP_s, r_s)$. The parameter $\eta$ is not important for our results and can be set equal to one. For the model described earlier on the lattice $\bZ^d$, results in the thermodynamic regime such as phase transitions are stated in terms of $\eta$. For the continuum model such results can be stated either in terms of $\eta$ or the intensity of the underlying Poisson process.

\section{Statements of main results}
\lab{Statements of main results}

Let $D_k = D_{k,s}$ be the number of vertices of degree $k$ in $G(\cP_s, r_s)$. We derive a scaling regime in which $E[D_k]$ converges under certain conditions on the parameters. In this regime we show that $D_k$ converges in distribution to a Poisson random variable under some additional conditions on the parameters. For fixed $\xi \in \bR$ consider the scaling $r_s \equiv r_s(\xi)$, defined by 

\begin{equation}
r_s(\xi)^d=\f{1}{c_0 s} \left(\log s +(k-1) \log\log s +\xi + \log\left(\f{\al\be}{k! d}\right) \right)
\lab{eq:scaling}
\end{equation}
where $c_0= \f{2\pi\al\be}{d(\al\be- d)} \eta^{\f{d}{\al}} \Gamma\left(1-\f{d}{\al} \right)$. For $\lam>0$, let  $ Po(\lam)$ be a Poisson random variable with mean $\lam$ and $\xrightarrow{d}$ denote convergence in distribution.
\begin{thm} 
Consider the random graph $G(\cP_s, r_s)$ with the connection function $g_s$ of the form~(\ref{eq:connection_function_s}). Suppose $\al > d$, $\beta > 1$ and the scaling parameter $r_s$ is as defined by (\ref{eq:scaling}). Then for any $k\geq 0$ we have
\begin{equation}
E[D_k]\to e^{-\xi} \mbox{ as } s \to \infty.
\lab{eq:exp_d_k}
\end{equation}
\lab{thm:Ex_poisson_conv}
\end{thm}
Note that we need $\min\{\al, \al \beta\} > d$ for the vertices to have finite degrees. This holds since we also require that the weights have a finite mean $(\beta > 1)$. Our next result shows that under some additional conditions, $D_k$ converges in distribution to a Poisson random variable. Note that the condition $\al \be > 2d$ in Theorem~\ref{thm:poisson_conv} is the condition required for the vertices of the graph to have degree with finite variance \cite{Deprez2018}.
\begin{thm} 
Consider the random graph $G(\cP_s, r_s)$ with the connection function $g_s$ of the form~(\ref{eq:connection_function_s}). Suppose $\al > d$, $\beta > 1$ and the scaling parameter $r_s$ is as defined by (\ref{eq:scaling}). If for any $k\geq 0$,  $\al\be > \max\{ 2d, (k+1)d, (2k+3)(\al-d)\}$, then,
\begin{equation}
D_k\xrightarrow{d} Po(e^{-\xi}) \mbox{ as } s \to \infty.
\lab{eq:poisson_D_k}
\end{equation}
\lab{thm:poisson_conv}
\end{thm}
Coming to connectivity, for random geometric graphs, that is, when the connection function is the indicator function on the unit ball, it is known (see 13.37, \citep{Penrose2003}) that the graph is connected with high probability in the absence of isolated nodes. Such a result would imply that under the conditions of Theorem~\ref{thm:poisson_conv} with $k=0$, $G(\cP_s, r_s)$ is connected with probability close to $\exp(-e^{-\xi})$ for large $s$. Though such a result is out of our reach at the moment, the next result provides a sufficient condition for the graph $G(\cP_s, \cdot)$ to be connected with high probability. To state the result we need some notation. The idea is to choose the scaling parameter in such a way that there is a one hop path (a two path) from each vertex to every one of its neighbours within a distance that guarantees connectivity in the usual random geometric graph with high probability. Let $E^W$ denote the expectation with respect to the weight $W$. $B(x,r)$ be the unit ball of radius $r$ centered at $x$ with respect to the Euclidean metric and the origin by $O$. For $\gamma >0$, let
\begin{equation}
\hat{r}_s(\gamma)^d := \f{\gamma \log s}{\kappa s}, \qquad \mbox{where} \qquad \kappa := \int_{B(O , 1)}  E^{W_z}\left[1-\exp\left(- \f{\eta W_z}{|z|^{\al}}\right)\right]\,dz.
\lab{eq:def_kappa}
\end{equation}
Note that $\kappa < \infty$ if $\beta > 1$ and $\al > d$. Let $\hat{g}_s(x,W_x; y,W_y)$ be the connection probability between two points $x,y$ with $r_s$ replaced by $\hat{r}_s(\gamma)(\equiv \hat{r}_s)$ in (\ref{eq:connection_function_s}).  Denote the random graph with the connection function $\hat{g}_s$ by $G(\cP_s,\hat{r}_s)$. Let $\th_d$ be the volume of the unit ball. Define the functions
\begin{equation}
T(\gamma):=   1-\exp\left\{ -\eta \left(1 + \left( \f{ \kappa }{\gamma \th_d} \right)^{\f{1}{d}}\right)^{-\al} \right\} \qquad \mbox{ and } \qquad Q(\gamma):= \gamma T(\gamma).
\lab{eq:def_TQ}
\end{equation}
Observe that $Q(0)=0, Q(\infty)=\infty$ and $Q'(\gamma)>0$ so that $\rho$ satisfying $Q(\rho) = 1$ is uniquely defined.

\begin{thm}
Let $Q$ be as defined in (\ref{eq:def_TQ}) and $\rho$ satisfy $Q(\rho) = 1$. Suppose $\al > d$, $\beta > 1$. Then for the sequence of graphs $G(\cP_s, \hat{r}_s(\gamma))$
\begin{equation}
P\left( G(\cP_s, \hat{r}_s(\gamma)) \mbox{ is connected}\right)\to 1,
\lab{eq:connectivity_eq}
\end{equation}
as $ s\to\infty$, for all $\gamma>\rho$.
%where $\rho$ is  defined as in~(\ref{eq:chi}).
\lab{thm:connectivity}
\end{thm}

\section{Proofs}

In all the proofs $c, c_0, c_1,c_2,\cdots$ and $C, C_0, C_1,C_2,\cdots$ will denote constants whose values will change from place to place. Let $d_1$ be the toroidal metric on $r_s^{-1}S$ defined as 
\begin{equation}
d_1(x,y) := \inf\limits_{z \in r_s^{-1} \bZ^d} ||x-y-z||, \qquad x,y \in r_s^{-1}S.
\lab{eqn:d1}
\end{equation}
In the proofs we often make the change of variable from $r_s^{-1}x$ to $x$ and $r_s^{-1}y$ to $y$ and hence shall refer to this as the standard change of variables. Such a change of variables transforms the connection function
$1 - \exp\left( - \f{\eta W_x W_y}{(d(x,y)/r_s)^\al} \right)$ to $1 - \exp\left( - \f{\eta W_x W_y}{(d_1(x,y))^\al} \right)$
%= 1 - \exp\left( - \f{\eta W_x W_y}{|y|^\al} \right).
%
which is independent of $s$ and we shall write as $\tilde{g}(x,W_x; y,W_y)$. 
%We denote the origin by $O$ and the ball of radius $r$ centered at $x$ by $B(x, r)$.
%
\subsection{Proof of Theorem~\ref{thm:Ex_poisson_conv} }
Fix $\xi \in \bR$ and let $r_s$ be as defined in (\ref{eq:scaling}). By the Campbell-Mecke formula we have
\begin{eqnarray}
E[D_k] & = & E\left[\sum_{X\in\cP_s}1_{\{\deg(X)=k \mbox{ in } G(\cP_s, r_s) \}}\right]=  s P^o\Big( \deg(O)=k \mbox{ in } G(\cP_s, r_s) \Big)\nn\\
&=&  s E^{W_0}\left[ \f{1}{k!}\left( s \int_{S} E^{W_x}\left[g_s \left(O, W_0; x, W_x\right)\right]\, dx\right)^k  \exp\left( -s \int_{S} E^{W_x}\left[g_s \left(O, W_0; x, W_x \right)\right]\, dx\right)  \right].
\lab{eq:translation_inv}
\end{eqnarray}
By the standard change of variables we obtain
\begin{eqnarray}
E[D_k] &=&  s E^{W_0}\left[ \f{1}{k!}\left( s r_s^d\int_{r_s^{-1} S} E^{W_x}\left[ \tilde{g}\left(W_0; x, W_x\right)\right]\, dx\right)^k  e^{ -s r_s^d \int_{r_s^{-1}S} E^{W_x}\left[ \tilde{g}\left(W_0; x, W_x\right)\right]\, dx}  \right],
\lab{eq:exp_scaled}
\end{eqnarray}
where we have written $E^{W_x}\left[\tilde{g} \left(W_0; x, W_x \right)\right]$ for $E^{W_x}\left[ \tilde{g}(O, W_0; x, W_x)\right]$. 
% 
%Using the above observation  we obtain
%
%\begin{equation}
%E[D_k] =  s E^{W_0}\left[ \f{1}{k!}  \left( -sr_s^d \int_{r_s^{-1} S} E^{W_x}\left[ g \left(W_0; x, W_x\right) \right]\, dx\right)^k  \exp\left( -s r_s^d\int_{r_s^{-1}S} E^{W_x}\left[ g \left(W_0; x, W_x\right)\right]\,dx\right)  \right].
%lab{eq:exp_scaled}
%\end{equation}
%
We shall compute upper and lower bounds for the expression on the right in (\ref{eq:exp_scaled}). Let $D_1:=B(O, \f{1}{2})$. We start with the following trivial inequalities.
\begin{equation}
\int_{r_s^{-1}D_1}E^{W_x}\left[ \tilde{g} \left(W_0; x, W_x\right) \right]\,dx
\leq \int_{r_s^{-1}S}E^{W_x}\left[ \tilde{g} \left(W_0; x, W_x\right) \right]\,dx
\leq \int_{\bR^d}E^{W_x}\left[ \tilde{g} \left(W_0; x, W_x\right) \right]\,dx.
\lab{eq:inner_outer_circle}
\end{equation}
Consider the upper bound in (\ref{eq:inner_outer_circle}). Using the probability density function of the weights and the fact that $d_1(O,x) = |x|$ we obtain
\[ \int_{\bR^d}E^{W_x}\left[ \tilde{g} \left(W_0; x, W_x\right) \right] \,dx  =  \int_{\bR^d}\int_{1}^{\infty}\left(1-\exp\left(-\f{\eta W_0 w}{|x|^\al}\right)\right) \be\, w^{-\be-1}\,dw\,dx \]
Switching to polar coordinates and applying the Fubini's theorem and then making the change of variables $t = r^{-\alpha}$ in the above equation yields
\begin{eqnarray}
\int_{\bR^d}E^{W_x}\left[ \tilde{g} \left(W_0; x, W_x\right) \right] \,dx  & = & \int_{0}^{\infty}\int_{0}^{2\pi}\int_{1}^{\infty}r^{d-1}\Big(1-\exp\Big(-\f{\eta W_0 w}{r^\al}\Big)\Big) \be\, w^{-\be-1} \,dw\,d\th\,dr\nn\\
& = &\f{2\pi \be}{\al} \int_{1}^{\infty}\int_{0}^{\infty}t^{-\f{d}{\al}-1}\Big(1-\exp\Big(-\eta W_0 w t\Big)\Big) \, w^{-\be-1}\,dt\,dw.
\lab{eq:outer_circle_1}
\end{eqnarray}
By a change of variable the inner integral in (\ref{eq:outer_circle_1}) can be evaluated to yield 
\begin{eqnarray}
\int_{0}^{\infty}t^{-\f{d}{\al}-1}\Big(1- e^{-\eta W_0 w t} \Big) \,dt 
& = & -\f{\al}{d} \left(1-e^{-\eta W_0 w t} \right) t^{-\f{d}{\al}} \bigg\vert_{0}^{\infty} +   \f{\al}{d} (\eta W_0 w) \int_{0}^{\infty} t^{-\f{d}{\al}}  e^{-\eta W_0 w t} \,dt \nn\\
& = &   \f{\al}{d}(\eta W_0 w)^{\f{d}{\al}} \int_{0}^{\infty}  u^{-\f{d}{\al}}  e^{- u} \,du =  \f{\al}{d}\left( \eta W_0 w \right)^{\f{d}{\al}} \Gamma\left(1-\f{d}{\al} \right),
\lab{eq:gamma_1}
\end{eqnarray}
where we  used the assumption that $\al > d$. 
Substituting from (\ref{eq:gamma_1}) in (\ref{eq:outer_circle_1}) we obtain 
\begin{eqnarray}
 \int_{\bR^d}E^{W_x}\left[ \tilde{g} \left(W_0; x, W_x\right) \right]\,dx  & = & \f{2\pi\be}{d}\left( \eta W_0\right)^{\f{d}{\al}} \Gamma\left(1-\f{d}{\al} \right) \int_{1}^{\infty}  w^{\f{d}{\al}-\be-1}=  c_0 W_0^{\f{d}{\al}},
\lab{eq:gamma_2}
\end{eqnarray}
where we have used the fact that $\al\be >d$ and set $c_0= \f{2\pi\al\be}{d(\al\be- d)} \eta^{\f{d}{\al}} \Gamma\left(1-\f{d}{\al} \right)$.
We now bound the integral on the left in (\ref{eq:inner_outer_circle}) from below. Proceeding as above and using (\ref{eq:gamma_2}) we obtain
\begin{equation}
\int_{r_s^{-1}D_1}E^{W_x}\left[ \tilde{g} \left(W_0; x, W_x\right) \right]\,dx =  c_0 W_0^{\f{d}{\al}} - 2 \pi \beta I_2, %\int_{r_s^{-1}D_1}\int_{1}^{\infty}\Big(1-\exp\Big(-\f{\eta W_0 w}{|x|^\al}\Big)\Big)\be\,w^{-\be-1}\,dw\,dx 
%\nn\\
 %
% & = &  \int_{0}^{\f{1}{2r_s}}\int_{0}^{2\pi}\int_{1}^{\infty}r^{d-1}\Big(1-\exp\Big(-\f{\eta W_0  w}{r^\al}\Big)\Big)\be\, w^{-\be-1}\,dw\,d\th\,dr  = 2\pi\be (I-I_2),
 %
\lab{eq:inner_circle_1}
\end{equation}
where by a use of the Fubini's theorem we can write
%
%$$I:= \int_{0}^{\infty}\int_{1}^{\infty}r^{d-1}\Big(1-\exp\Big(-\f{\eta W_0 w}{r^\al}\Big)\Big)\, w^{-\be-1}\,dw\,dr \qquad \mbox{and} \qquad
%$$I_2:=\int_{\f{1}{2r_s}}^{\infty}\int_{1}^{\infty}r^{d-1}\Big(1-\exp\Big(-\f{\eta W_0  w}{r^\al}\Big)\Big)\, w^{-\be-1}\,dw\,dr.$$ 
%
\begin{equation}
I_2  = \int_{1}^{\infty}\int_{\f{1}{2r_s}}^{\infty}\Big[1-\exp\Big(-\f{\eta W_0 w}{r^\al}\Big)\Big]r^{d-1} \, w^{-\be-1}\,dr\,dw.
\lab{eq:I_2}
\end{equation}
To compute the inner integral on the right hand side in (\ref{eq:I_2}) we make the change of variable $t = r^{-\alpha}$ and use the assumption that $\al > d$.
\begin{eqnarray}
%
% \lefteqn{ \int_{\f{1}{2r_s}}^{\infty}\left[1-\exp\Big(-\f{\eta W_0 w}{r^\al}\Big)\right]r^{d-1}\,dr  
 \f{1}{\al}\int_{0}^{(2r_s)^\al} \left(1-e^{-\eta W_0 w t}\right) t^{-\f{d}{\al} - 1}\,dt
& = & -\f{1}{d} \left(1-e^{-\eta W_0 w t} \right) t^{-\f{d}{\al}} \bigg\vert_{0}^{(2r_s)^\al} + \f{1}{d}(\eta W_0 w) \int_{0}^{(2r_s)^\al} t^{-\f{d}{\al}}  e^{-\eta W_0 w t} \,dt  \nn \\
%
%& = & - \f{1}{d (2r_s)^d}\left[ 1-\exp(-\eta W_0 w (2r_s)^\al) \right] +  \f{1}{d}(\eta W_0 w) \int_{0}^{(2r_s)^\al}  t^{-\f{d}{\al}}  \exp(- \eta W_0 w t) \,dt\nn\\
%
& \leq &   \f{1}{d}(\eta W_0 w) \int_{0}^{(2r_s)^\al}  t^{-\f{d}{\al}}  \,dt
\; = \;   c_1 W_0 w r_s^{\al-d},
\lab{eq:I_2_1}
\end{eqnarray}
where $c_1=\f{\al\eta 2^{\al-d}}{d(\al-d)}$. From (\ref{eq:I_2}), (\ref{eq:I_2_1}) and the fact that $\beta > 1$ we obtain 
\begin{eqnarray}
I_2 \leq \int_{1}^{\infty}c_1 \be W_0 r_s^{\al-d} w^{-\be}= c_2 W_0 r_s^{\al-d},
\lab{eq:I_2_upperbound}
\end{eqnarray}
where $c_2= \f{\be c_1}{\be-1}$. From (\ref{eq:inner_circle_1}), (\ref{eq:I_2_upperbound}) and the fact that $\tilde{g} > 0$ we obtain
\begin{equation}
\int_{r_s^{-1}D_1} E^{W_x}\left[ \tilde{g} \left(W_0; x, W_x\right) \right]\,dx\geq 0 \vee \left(  c_0 W_0^{\f{d}{\al}}- c_3 W_0 r_s^{\al-d}\right),
\lab{eq:lower_bound_inner_circle}
\end{equation}
where $c_3= 2\pi \beta c_2$.  For $w\geq 1$, let $\Lam_s(w):= 0 \vee \left(  c_0 w^{\f{d}{\al}}- c_3 w r_s^{\al-d}\right)$. Substituting from (\ref{eq:gamma_2}) and (\ref{eq:lower_bound_inner_circle}) in (\ref{eq:inner_outer_circle}) we obtain 
\begin{equation}
%0 \vee \left(  c_0 W_0^{\f{d}{\al}}- c_3 W_0 r_s^{\al-d}\right)
%
\Lam_s(W_0) \leq \int_{r_s^{-1}S}E^{W_x}\left[ \tilde{g} \left(W_0; x, W_x\right) \right]\,dx\leq c_0 W_0^{\f{d}{\al}}.
\lab{eq:in_out_circle}
\end{equation}
It follows from (\ref{eq:exp_scaled}) and (\ref{eq:in_out_circle}) that
\begin{equation}
 \f{s}{k!} E^{W_0} \left[    \left( s r_s^d \Lam_s(W_0)\right)^k  e^{- c_0 s r_s^d  W_0^{\f{d}{\al}}} \right]  \leq
E[D_k]   \leq  \f{s}{k!}  E^{W_0} \left[  \left( c_0 s r_s^d  W_0^{\f{d}{\al}} \right)^k   e^{- s r_s^d\Lam_s(W_0)} \right].
\lab{eq:expectation_bound_1}
\end{equation}
% % % % %
%\begin{eqnarray}
%E[D_k] & \geq & \f{s}{k!} E^{W_0} \left[    \left( s r_s^d \Lam_s(W_0)\right)^k  e^{- c_0 s r_s^d  W_0^{\f{d}{\al}}} \right],  \nn \\
% E[D_k]  & \leq & \f{s}{k!}  E^{W_0} \left[  \left( c_0 s r_s^d  W_0^{\f{d}{\al}} \right)^k   e^{- s r_s^d\Lam_s(W_0)} \right].
%  
%\lab{eq:expectation_bound_1}
%\end{eqnarray}
% % % % % %
The next step is to prove that both the upper and lower bounds in (\ref{eq:expectation_bound_1}) converge to $e^{-\xi}$ as $s\to\infty$. For later use we shall state the convergence of the lower and upper bounds in somewhat greater generality as separate lemmas.

Theorem~\ref{thm:Ex_poisson_conv} now follows from (\ref{eq:expectation_bound_1}) and Lemmas~\ref{lem:equal_limit_L}--\ref{lem:equal_limit_U} stated below.  For the case $k=0$ we use the first assertion in Lemma~\ref{lem:equal_limit_L} and Lemma~\ref{lem:equal_limit_UB} with $j=1$. For $k \geq 1$ we use the second assertion in Lemma~\ref{lem:equal_limit_L} and Lemma~\ref{lem:equal_limit_U} with $j=1$ and $m=1$. \qed
\begin{lem}
\begin{enumerate}[(i)]
\item \lab{equal_limit_L1} Let $r_s$ be as defined in (\ref{eq:scaling}) with $k=0$.
\begin{equation} 
\lim\limits_{s \to \infty} s E^{W_0} \left[   \exp{\left(- c_0 s r_s^d  W_0^{\f{d}{\al}}\right)} \right] = e^{-\xi}.
\lab{eq:equal_limit_L1}
\end{equation}
\item \lab{equal_limit_L2} Let $r_s$ be as defined in (\ref{eq:scaling}). Suppose $\al>d$ and $k\geq 1$, then
\begin{equation} 
\lim\limits_{s \to \infty} \f{s}{k!} E^{W_0} \left[    \left( s r_s^d \Lam_s(W_0)\right)^k  \exp{\left(- c_0 s r_s^d  W_0^{\f{d}{\al}}\right)} \right] = e^{-\xi}.
\lab{eq:equal_limit_L2}
\end{equation}
\end{enumerate}
\lab{lem:equal_limit_L}
\end{lem}
\begin{lem}
Let $r_s$ be as defined in (\ref{eq:scaling}) with $k=0$. Suppose $j\geq 1$, $\al>d$ and $\al\be>jd$, then as $s\to \infty$
%for all suffuciently large $s$
%
\begin{equation}
s^j E^{W_0} \left[  e^{ - j s r_s^d \Lam_s(W_0) } \right] = \f{1}{j} \left(\f{\al\be}{d}\right)^{-j+1} e^{-j\xi} (\log s)^{j-1} +o(1).
\lab{eq:equal_limit_UB}
\end{equation}
\lab{lem:equal_limit_UB}
\end{lem}
\begin{lem}
Let $r_s$ be as defined in (\ref{eq:scaling}) with $k\geq 1$. Suppose $\al> d $, $j\geq 1$  and $m\geq 1$, then as $s\to \infty$
%for all suffuciently large $s$
%
\begin{equation}
\left(\f{s}{k!}\right)^j E^{W_0} \left[  \left( c_0 s r_s^d  W_0^{\f{d}{\al}} \right)^{jm} e^{ - j s r_s^d \Lam_s(W_0) }\right] = \f{1}{j} \left(\f{\al\be}{d}\right)^{-j+1} e^{-j\xi} (\log s)^{j-1} +o(1).
\lab{eq:equal_limit_U_3}
\end{equation}
\lab{lem:equal_limit_U}
\end{lem}
The following observations will be invoked several times in the proofs. 
\begin{equation}
c_0 s r_s^d = \log s + (k-1)\log\log s +\xi + \log \left(\f{\al\be}{k! d }\right), \qquad \left(c_0 s r_s^d \right)'= \f{d}{ds}\left(c_0 s r_s^d \right) = s^{-1}(1 + o(1)).
\lab{eq:useful_observation}
\end{equation}

\textbf{Proof of Lemma~\ref{lem:equal_limit_L}(\ref{equal_limit_L1}).} Using the density of $W_0$ we obtain
\begin{eqnarray}
s E^{W_0} \left[  \exp\left( - c_0 s r_s^d  W_0^{\f{d}{\al}} \right) \right]  
&=& s \int_{1}^{\infty} \exp\left( - c_0 s r_s^d  w^{\f{d}{\al}} \right) \be w^{-\be-1}\, dw \nn\\
& = &   \f{\al\be}{d}s\, \left(c_0 s r_s^d\right) ^{ \f{\al\be}{d}} \int_{c_0 s r_s^d}^{\infty} e^{ - t} t^{-\f{\al\be}{d}-1}  \, dt 
=  \f{\al\be}{d} \, \f{\int_{c_0 s r_s^d}^{\infty} e^{ - t} t^{-\f{\al\be}{d}-1}  \, dt}{s^{-1}\left(c_0 s r_s^d\right) ^{ -\f{\al\be}{d}}}. 
\lab{eq:equal_limit_3a}
\end{eqnarray}
Since $sr_s^d \to \infty$ as $s \to \infty$ we can apply to L'Hospital's rule to conclude that the limit of the last expression in (\ref{eq:equal_limit_3a}) equals the limit of the ratio
%
%\begin{eqnarray}
\[ \f{\al\be}{d} \, \f{-  e^{ - c_0 s r_s^d} (c_0 s r_s^d)^{-\f{\al\be}{d}-1} \left(c_0 sr_s^d \right)'}{ -s^{-2}(c_0 s r_s^d)^{-\f{\al\be}{d}} - \f{\al\be}{d}  s^{-1}  (c_0 s r_s^d)^{-\f{\al\be}{d}-1} \left(c_0 sr_s^d \right)' }  
= \f{\al\be}{d} \, \f{ s^2 e^{ - c_0 s r_s^d} \left(c_0 sr_s^d \right)'}{ (c_0 s r_s^d) + s \f{\al\be}{d} \left(c_0 sr_s^d \right)'}. \]
%
%\lab{eq:equal_limit_3}
%\end{eqnarray}
%
Lemma~\ref{lem:equal_limit_L}(\ref{equal_limit_L1}) now follows from (\ref{eq:useful_observation}) with $k=0$.
%by substituting $c_0 s r_s^d = \log s+\xi - \log\log s + \log \left(\f{\al\be}{d}\right)$ and observing that $\left(c_0 sr_s^d \right)' = s^{-1}(1 + o(1))$. \qed
%

\textbf{Proof of Lemma~\ref{lem:equal_limit_L}(\ref{equal_limit_L2}).} Substituting for the density of $W_0,$ the expression in the limits on the left of (\ref{eq:equal_limit_L2}) equals
%Using the density of $W_0$ we obtain
%
\begin{eqnarray}
%
%\lefteqn{ \f{s}{k!} E^{W_0} \left[   \left( s r_s^d \left(0 \vee \left(  c_0 W_0^{\f{d}{\al}}- c_3 W_0 r_s^{\al-d}\right) \right)\right)^k  \exp\left( - c_0 s r_s^d  W_0^{\f{d}{\al}} \right) \right]  }\nn\\
%
%&=& 
\lefteqn{ \f{s}{k!} \int_{1}^{\infty}  \left( s r_s^d \left(0 \vee \left(  c_0 w^{\f{d}{\al}}- c_3 w r_s^{\al-d}\right) \right)\right)^k  \exp\left( - c_0 s r_s^d  w^{\f{d}{\al}} \right) \be w^{-\be-1}\, dw} \nn\\
&=& \f{s}{k!} \int_{1}^{C r_s^{-\al}}  \left( c_0 s r_s^d w^{\f{d}{\al}}- c_3 s w r_s^{\al}\right)^k  \exp\left( - c_0 s r_s^d  w^{\f{d}{\al}} \right) \be w^{-\be-1}\, dw,
\lab{eq:equal_limit_3ab}
\end{eqnarray}
where $C= \left(\f{c_0}{c_3}\right)^{\f{\al}{\al-d}}$. 
By the change of variable $t= c_0 s r_s^d w^{\f{d}{\al}}$ the right hand side in (\ref{eq:equal_limit_3ab}) equals
\begin{equation}
 \f{\al\be}{k! d}s\, \left(c_0 s r_s^d\right) ^{ \f{\al\be}{d}} \int_{c_0 s r_s^d}^{c_4 s} \left( t - c_3 c_0^{-\f{\al}{d}}\f{t^{\f{\al}{d}}}{s^{\f{\al}{d}-1}} \right)^k e^{ - t} t^{-\f{\al\be}{d}-1}  \, dt
=  \f{\al\be}{k! d} \, \f{\int_{c_0 s r_s^d}^{c_4 s} \left( t - c_5\f{t^{\f{\al}{d}}}{s^{\f{\al}{d}-1}} \right)^k  e^{ - t} t^{-\f{\al\be}{d}-1}  \, dt}{s^{-1}\left(c_0 s r_s^d\right) ^{ -\f{\al\be}{d}}}, 
\lab{eq:equal_limit_3ac}
\end{equation}
where $c_4 = c_0^{\f{\al}{\al-d}}c_3^{-\f{d}{\al-d}}$ and $c_5=c_3 c_0^{-\f{\al}{d}}$. 
Again we can apply L'Hospital's rule to conclude that the limit of the ratio on the right in (\ref{eq:equal_limit_3ac}) equals the limit of
\begin{equation}
\f{\al\be}{k!d} \left[ \f{  \left( c_4 s - c_5 c_4^{\f{\al}{d}} s  \right)^k e^{ - c_4 s } (c_4 s )^{-\f{\al\be}{d}-1} c_4 }{-s^{-2}(c_0 s r_s^d)^{-\f{\al\be}{d}} -  \f{\al\be}{d} s^{-1}  (c_0 s r_s^d)^{-\f{\al\be}{d}-1} \left(c_0 sr_s^d \right)'} 
 - \f{ \left( c_0 s r_s^d - c_3 s r_s^{\al} \right)^k e^{ - c_0 s r_s^d} (c_0 s r_s^d)^{-\f{\al\be}{d}-1} \left(c_0 sr_s^d \right)'}{ -s^{-2}(c_0 s r_s^d)^{-\f{\al\be}{d}} - \f{\al\be}{d}  s^{-1}  (c_0 s r_s^d)^{-\f{\al\be}{d}-1} \left(c_0 sr_s^d \right)' } \right]
\lab{eq:equal_limit_3}
\end{equation}
The first term in (\ref{eq:equal_limit_3}) is zero since $c_4 - c_5 c_4^{\f{\al}{d}}=0$. The second term (\ref{eq:equal_limit_3}) simplifies to 
%
%\begin{eqnarray}
%
\[ \f{\al\be}{k! d} \, \f{ s^2  \left( c_0 s r_s^d - c_3 s r_s^{\al} \right)^k e^{ - c_0 s r_s^d} \left(c_0 sr_s^d \right)'}{ (c_0 s r_s^d) + s \f{\al\be}{d} \left(c_0 sr_s^d \right)'}
= \f{\al\be}{k! d} \, \f{ s^2  \left( c_0 s r_s^d  \right)^k e^{ - c_0 s r_s^d} \left(c_0 sr_s^d \right)'}{ (c_0 s r_s^d) + s \f{\al\be}{d} \left(c_0 sr_s^d \right)'} \left( 1 -\f{ c_3 }{ c_0 } r_s^{\al-d} \right)^k. \]
%
%\lab{eq:equal_limit_4}
%\end{eqnarray}
%
Lemma~\ref{lem:equal_limit_L}(\ref{equal_limit_L2}) now follows from (\ref{eq:useful_observation})
%by substituting $c_0 s r_s^d = \log s + (k-1)\log\log s +\xi + \log \left(\f{\al\be}{k! d }\right)$ 
and observing that 
%$\left(c_0 sr_s^d \right)' = s^{-1}(1 + o(1)),$ 
$\al>d$ and $r_s \to 0$. \qed

\textbf{Proof of Lemma~\ref{lem:equal_limit_UB}.}
Observe that $\{c_0 W_0^{\f{d}{\al}}- c_3 W_0 r_s^{\al-d} \geq 0\}=\{W_0 \leq C r_s^{-\al}\}$ where $C= \left(\f{c_0}{c_3}\right)^{\f{\al}{\al-d}}$.  The left hand side of (\ref{eq:equal_limit_UB}) equals
\begin{eqnarray}
%\lefteqn{s^j  E^{W_0} \left[  \exp\left(- j s r_s^d\max{\left\{0, c_0 W_0^{\f{d}{\al}}- c_3 W_0 r_s^{\al-d}\right\}} \right) \right] }\nn\\
%
%& = &  
\lefteqn{s^j E^{W_0} \left[  \exp\left(- j c_0 s r_s^d W_0^{\f{d}{\al}}+ j c_3 W_0 s r_s^{\al} \right) ; W_0 \leq C r_s^{-\al}\right] +   s ^j P\left(W_0 > C r_s^{-\al}\right)}\nn\\
& = &  s^j \int_{1}^{C r_s^{-\al}} \exp \left( - j c_0 s r_s^d w^{\f{d}{\al}} + j c_3 sr_s^\al w\right)\be w^{-\be-1}\, dw +  s^j C^{-\be } r_s^{\al\be}.
\lab{eq:equal_limit_5a}
\end{eqnarray}
The second term in (\ref{eq:equal_limit_5a}) converges to zero by (\ref{eq:scaling}) since $\al \beta > jd$. By changing the variable $t= j c_0 s r_s^d w^{\f{d}{\al}}$ in the first term in (\ref{eq:equal_limit_5a}) we obtain 
\begin{equation}
\f{\al\be}{d} s^j \, \left(j c_0 s r_s^d\right)^{ \f{\al\be}{d}} \int_{j c_0 s r_s^d}^{c_4s} e^{- t+c_5 t^{\f{\al}{d}} s^{-\f{\al}{d}+1}} t^{-\f{\al\be}{d}-1}  \, dt  
=  \f{\al\be}{d}   \f{\int_{j c_0 s r_s^d}^{c_4s} e^{ - t+c_5 t^{\f{\al}{d}} s^{-\f{\al}{d}+1}} t^{-\f{\al\be}{d}-1}  \, dt}{s^{-j}\left(j c_0 s r_s^d\right) ^{ -\f{\al\be}{d}}},
\lab{eq:equal_limit_5_1}
\end{equation}
where $c_4= j c_0 C^{\f{d}{\al}}$ and $c_5=j c_3 (j c_0)^{-\f{\al}{d}}$. By the L'Hospital's rule the limit of the expression on the right in (\ref{eq:equal_limit_5_1}) equals the limit of
\begin{eqnarray}
 \lefteqn{\f{\al\be}{d}   \f{e^{-c_4 s+c_5 (c_4 s)^{\f{\al}{d}} s^{-\f{\al}{d}+1}} (c_4 s)^{-\f{\al\be}{d}-1} c_4    -  e^{ - j c_0 s r_s^d } (j c_0 s r_s^d)^{-\f{\al\be}{d}-1} \left( j c_0 s r_s^d \right)'  e^{c_5 (jc_0)^{\f{\al}{d}} s r_s^\al}  }{ -j s^{-j-1}(j c_0 s r_s^d)^{-\f{\al\be}{d}} - \f{\al\be}{d}  s^{-j}  (j c_0 s r_s^d)^{-\f{\al\be}{d}-1} \left( j c_0 sr_s^d \right)' }}\nn\\
& = &   C_1  \f{e^{-c_4 s+c_5 c_4^{\f{\al}{d}} s} s^{-\f{\al\be}{d}+j} ( j c_0 s r_s^d)^{\f{\al\be}{d}+1} }{  j ( j c_0 s r_s^d) + s \f{\al\be}{d} \left( j c_0 sr_s^d \right)' }
 +  \f{\al\be}{d}    \f{  s^{j+1} e^{ - j c_0 s r_s^d } \left(j c_0 s r_s^d \right)'}{ j ( j c_0 s r_s^d) + s \f{\al\be}{d} \left( j c_0 sr_s^d \right)'} e^{c_5 (jc_0)^{\f{\al}{d}} s r_s^\al},
\lab{eq:equal_limit_6_1}
\end{eqnarray}
where $C_1$ is a constant. Since $c_4  - c_5 c_4^{\f{\al}{d}}=0$, the first term on the right in (\ref{eq:equal_limit_6_1}) simplifies to
%
%\begin{eqnarray}
%
\[ C_1  \f{ s^{-\f{\al\be}{d}+j} ( j c_0 s r_s^d)^{\f{\al\be}{d}+1} }{ j ( j c_0 s r_s^d) + s \f{\al\be}{d} \left( j c_0 sr_s^d \right)' }\to 0 \mbox{ as } s \to \infty, \]
%
%\lab{eq:}
%\end{eqnarray}
%
by (\ref{eq:useful_observation}) and the assumption that $\al\be > jd$. 
Lemma~\ref{lem:equal_limit_UB} now follows by using (\ref{eq:useful_observation}) in the second term on the right in (\ref{eq:equal_limit_6_1}) and the assumption that $\al>d$. \qed

\textbf{Proof of Lemma~\ref{lem:equal_limit_U}.}
By the observation at the beginning of the proof of Lemma~\ref{lem:equal_limit_UB}, the left hand side of (\ref{eq:equal_limit_U_3}) equals
\begin{equation}
\left( \f{s}{k!}\right)^j \left[ \int_{1}^{C r_s^{-\al}} \left( c_0 s r_s^d  w^{\f{d}{\al}} \right)^{jm} e^{ - j c_0 s r_s^d w^{\f{d}{\al}} + j c_3 sr_s^\al w} \be w^{-\be-1}\, dw + \int_{C r_s^{-\al}}^{\infty} \left( c_0 s r_s^d  w^{\f{d}{\al}} \right)^{jm} \be w^{-\be-1}\, dw \right] .
\lab{eq:equal_limit_7}
\end{equation}
%\end{eqnarray}
%
By changing the variable $t= c_0 s r_s^d w^{\f{d}{\al}}$ in the second term in (\ref{eq:equal_limit_7}) we obtain 
\begin{eqnarray}
 \f{\al\be}{\left(k!\right)^j d} \; s^j \, \left( c_0 s r_s^d\right)^{\f{\al\be}{d}} \int_{c_0 C^{\f{d}{\al}} s }^{\infty} t^{jm- \f{\al\be}{d}-1}\, dt
= \f{\al\be}{\left(k!\right)^j  d} \f{ \int_{c_4 s}^{\infty} t^{jm- \f{\al\be}{d}-1}\, dt}{s^{-j} \, \left( c_0 s r_s^d\right)^{-\f{\al\be}{d}}},
\lab{eq:second_term_j}
\end{eqnarray}
where $c_4=  c_0 C^{\f{d}{\al}}$. To apply the L'Hospital's rule differentiate the numerator and denominator of the expression on the right in (\ref{eq:second_term_j}) to obtain
\begin{equation}
\f{\al\be}{\left(k!\right)^j  d}   \f{- (c_4 s)^{jm-\f{\al\be}{d}-1} c_4 }{ -j s^{-j-1}( c_0 s r_s^d)^{-\f{\al\be}{d}} - \f{\al\be}{d}  s^{-j}  ( c_0 s r_s^d)^{-\f{\al\be}{d}-1} \left( c_0 sr_s^d \right)' }=
c_4^{jm-\f{\al\be}{d}}  \f{\al\be}{\left(k!\right)^j  d}  \f{ s^{jm-\f{\al\be}{d}+j} ( c_0 s r_s^d)^{\f{\al\be}{d}+1} }{  j (c_0 s r_s^d) + s \f{\al\be}{d} \left(  c_0 sr_s^d \right)' }.
\lab{eq:second_term_j1}
\end{equation}
By making the same change of variable $t= c_0 s r_s^d w^{\f{d}{\al}}$ in the first term in (\ref{eq:equal_limit_7}) we obtain 
\begin{equation}
\f{\al\be}{\left(k!\right)^j  d} s^j \, \left( c_0 s r_s^d\right)^{ \f{\al\be}{d}} \int_{c_0 s r_s^d}^{c_4s} e^{ - j t+c_5 j t^{\f{\al}{d}} s^{-\f{\al}{d}+1}} t^{jm-\f{\al\be}{d}-1}  \, dt  
=   \f{\al\be}{\left(k!\right)^j  d}   \f{\int_{ c_0 s r_s^d}^{c_4s} e^{ - jt+c_5 j t^{\f{\al}{d}} s^{-\f{\al}{d}+1}} t^{jm-\f{\al\be}{d}-1}  \, dt}{s^{-j}\left( c_0 s r_s^d\right) ^{ -\f{\al\be}{d}}},
\lab{eq:equal_limit_7_1}
\end{equation}
where $c_5=c_3 c_0^{-\f{\al}{d}}$. Differentiating the numerator and denominator of the expression on the right in (\ref{eq:equal_limit_7_1}) we obtain 
\[
 \f{\al\be}{\left(k!\right)^j  d}  \f{e^{-c_4 j s+c_5 j (c_4 s)^{\f{\al}{d}} s^{-\f{\al}{d}+1}} (c_4 s)^{jm-\f{\al\be}{d}-1} c_4 -  e^{ - j c_0 s r_s^d  + j c_5 c_0^{\f{\al}{d}} s r_s^\al  } ( c_0 s r_s^d)^{jm-\f{\al\be}{d}-1} \left(  c_0 s r_s^d \right)' }{ -j s^{-j-1}( c_0 s r_s^d)^{-\f{\al\be}{d}} - \f{\al\be}{d}  s^{-j}  ( c_0 s r_s^d)^{-\f{\al\be}{d}-1} \left( c_0 sr_s^d \right)' }
\]
\begin{eqnarray}
\quad  & = &   \f{\al\be}{\left(k!\right)^j  d} \left[ - c_4^{jm-\f{\al\be}{d}}   \f{e^{-c_4 j s+c_5 c_4^{\f{\al}{d}} j s} s^{jm-\f{\al\be}{d}+j} ( c_0 s r_s^d)^{\f{\al\be}{d}+1} }{  j (c_0 s r_s^d) + s \f{\al\be}{d} \left(  c_0 sr_s^d \right)' }  +   \f{  s^{j+1} e^{ - j c_0 s r_s^d } ( c_0 s r_s^d)^{jm}\left( c_0 s r_s^d \right)'}{ j ( c_0 s r_s^d) + s \f{\al\be}{d} \left( c_0 sr_s^d \right)'} e^{j c_5 c_0^{\f{\al}{d}} s r_s^\al} \right]\nn\\
& = &   \f{\al\be}{\left(k!\right)^j  d} \left[ - c_4^{jm-\f{\al\be}{d}}   \f{ s^{jm-\f{\al\be}{d}+j} ( c_0 s r_s^d)^{\f{\al\be}{d}+1} }{ j ( c_0 s r_s^d) + s \f{\al\be}{d} \left( c_0 sr_s^d \right)' }
+     \f{  s^{j+1} e^{- j c_0 s r_s^d  } ( c_0 s r_s^d)^{jm}\left( c_0 s r_s^d \right)'}{ j ( c_0 s r_s^d) + s \f{\al\be}{d} \left( c_0 sr_s^d \right)'} e^{ j c_5 c_0^{\f{\al}{d}} s r_s^\al} \right],
\lab{eq:equal_limit_6}
\end{eqnarray}
since $c_4  - c_5 c_4^{\f{\al}{d}}=0$. Adding (\ref{eq:second_term_j1}) and (\ref{eq:equal_limit_6}) and observing that the expression on the right in (\ref{eq:second_term_j1}) cancels with the first term on the right in (\ref{eq:equal_limit_6}) we obtain by the L'Hospital's rule that the limit of the expression on the left in (\ref{eq:equal_limit_U_3}) equals the limit of 
\[ \f{\al\be}{\left( k! \right)^j d} \; \f{  s^{j+1} \exp\left( - j c_0 s r_s^d  \right) ( c_0 s r_s^d)^{jm}\left( c_0 s r_s^d \right)'}{ j ( c_0 s r_s^d) + s \f{\al\be}{d} \left( c_0 s r_s^d \right)'} \exp\left( j c_5 c_0^{\f{\al}{d}} s r_s^\al \right). \]
%
%\begin{eqnarray}
%\lefteqn{ \lim\limits_{s\to \infty} \left( \f{s}{k!}\right)^j  E^{W_0} \left[  \left( c_0 s r_s^d  W_0^{\f{d}{\al}} \right)^{jk} \exp\left(- j s r_s^d\left(0 \vee \left(  c_0 W_0^{\f{d}{\al}}- c_3 W_0 r_s^{\al-d}\right) \right)  \right) \right] }\nn\\
%
%& = & \lim\limits_{s\to \infty} \f{\al\be}{\left( k! \right)^j d}     \f{  s^{j+1} \exp\left( - j c_0 s r_s^d  \right) ( c_0 s r_s^d)^{jk}\left( c_0 s r_s^d \right)'}{ j ( c_0 s r_s^d) + s \f{\al\be}{d} \left( c_0 s r_s^d \right)'} \exp{( j c_5 c_0^{\f{\al}{d}} s r_s^\al)}.\nn
%
%\end{eqnarray}
%
Lemma~\ref{lem:equal_limit_U} now follows by using (\ref{eq:useful_observation}) and the assumption that $\al>d$. \qed

\subsection{Stein's Method for Poisson Convergence}
To prove the Poisson convergence results stated in Theorem~\ref{thm:poisson_conv} we adapt Theorem 3.1 from \citep{Penrose2018}. This result which is based on the Stein's method, is modified to account for the inhomogeneity arising from the random weights associated with the vertices of the graph.

Let $\tilde{\cP}_s$ be  a $W$-marking of $\cP_s$. That is, suppose $\cP_s = \sum_{i=1}^{N_s} \del_{X_n}$ where $N_s$ is a Poisson random variable with mean $s$ and $X_1, X_2, \ldots$ is a sequence of independent uniformly distributed random variables taking values in $S$. Then $\tilde{\cP}_s = \sum_{i=1}^{N_s} \del_{(X_n, W_n)}$ where $W_1, W_2, \ldots$ is a sequence of independent random variables with probability distribution satisfying  $P(W>w) = w^{- \be}1_{[1,\infty)}(w)$ for some $\be > 0$ and independent of the random variable $N_s$ and the sequence $\{X_n\}_{n \geq 1}$. We write
\begin{equation}
D_k = D_k(\tilde{\cP}_s) = \sum_{y \in \cP_s } f(y,W_y, \tilde{\cP}_s \setminus \{(y,W_y)\}),
\label{eq:D_0}
\end{equation}
where $f(y,W_y, \tilde{\cP}_s \setminus \{(y,W_y)\})$ equals one if the degree of $y$ equals $k$ in the graph $G(\cP_s\cup \{y\}, r_s)$ and zero otherwise. We shall abbreviate $f(y,W_y, \tilde{\cP}_s)$ to $f(y,\cP_s)$. Set $\tilde{p}_s(x, W_x)= E\left[ f(x, \cP_s) \big\vert  W_x \right]$, $x \in S$. Since the underlying point process is homogeneous and the metric is toroidal, $\tilde{p}$ depends on $s, W_x$ and not on $x$. Hence we shall write $\tilde{p}_s(W_x)$ instead of $\tilde{p}_s(x, W_x)$. By the Campbell-Mecke formula $\nu := E[D_k]$ satisfies
\begin{equation}
\nu =  s  \int_S E^{W_x}[\tilde{p}_s(x,W_x)] \, dx = s \, E^{W_o}[\tilde{p}_s(W_o)].
\lab{eq:def_nu}
\end{equation}
Let $d_{TV}$ denote the total variation distance, $Z_\nu$ be a Poisson random variable with mean $\nu$ and $F_{D_k}, F_{Z_\nu}$ denote the distribution functions of $D_k, Z_{\nu}$ respectively. For any function $\phi : \bN \cup \{0\} \to \bR$, let $\Delta \phi(i) = \phi(i+1) - \phi(i)$ and $||\cdot||_{\infty}$ denote the $\sup$ norm.
\begin{thm}
Suppose that for almost every $x \in S$ we can find a random variable $V_x = V_x(W_x)$ coupled with $D_k$ such that conditional on $W_x$
%
%\begin{enumerate}[(i)]
%
%\item $U_x  \bigg\vert W_x \stackrel{d}{=} D_k\left(\{(x,W_x)\} \cup \tilde{\eta}_s \right)$,
%
%\item $
\[ 1+V_x \stackrel{d}{=} D_k\left(\{(x,W_x)\} \cup \tilde{\cP}_s \right) \big\vert \left\{ f \left( x, \cP_s  \right)=1 \right\}, W_x. \]
%
%\item $E\left[ |U_x-V_x|  \big\vert W_x \right] \leq q(x, W_x)$.
%
%\end{enumerate}
%
Then
\begin{equation}
d_{TV}\left( F_{D_k}, F_{Z_{\nu}} \right)  \leq   (1 \wedge \nu^{-1}) s \int_S E^{W_x} \left[ E\left[ |D_k-V_x|  \big\vert W_x \right] \, \tilde{p}_s(W_x) \right] \, dx.
\lab{eq:stein_bound}
\end{equation}
\lab{thm:SteinsMethod_PC}
\end{thm}
\textbf{Proof of Theorem~\ref{thm:SteinsMethod_PC}.}
%First we estimate an upper bound for $\big\vert E[ \nu \phi(J+1) - J \phi(J) ] \big\vert$ for 
Let $\phi : \bN \cup \{0\} \to \bR$ be a bounded function. By using the definition of $D_k$ and the Campbell-Mecke formula we have
\begin{eqnarray}
E[ D_k \phi(D_k) ] & = & E  \left[ \sum_{x \in \cP_s } f(x, \cP_s ) \phi(D_k(\tilde{\cP}_s)) \right]\nn\\
& = &  s\int_S E^{W_x}\left[ E \left[  f( x, \cP_s ) \phi\left(D_k(\{(x,W_x)\} \cup \tilde{\cP}_s)\right) \big\vert W_x \right] \right] \, dx \nn\\
& = & s \int_S E^{W_x} \left[ E \left[ \phi\left(D_k(\{(x,W_x)\} \cup \tilde{\cP}_s ) \right)  \big\vert   \{ f( x,  \eta_s )=1\} , W_x \right] \tilde{p}_s(W_x) \right] \, dx \nn \\
& = & s \int_S E^{W_x} \left[ E\left[  \phi(V_x+1)  \big\vert W_x \right]  \tilde{p}_s(W_x) \right] \, dx.
\lab{eq:Wg(W)}
\end{eqnarray}
From (\ref{eq:def_nu}) we have
\begin{eqnarray}
E\left[ \nu \phi(D_k+1) \right] & = &  s\int_S E^{W_x} \left[ \phi(D_k+1)  \tilde{p}_s(W_x) \right]\, dx.
\lab{eq:alpha_g(W+1)}
\end{eqnarray}
Using~(\ref{eq:Wg(W)}) and~(\ref{eq:alpha_g(W+1)}) we obtain
\begin{eqnarray}
\big\vert E\left[ \nu \phi(D_k+1) - D_k \phi(D_k) \right] \big\vert & \leq & 
 s  \int_S E^{W_x} \left[  \big\vert E\left[  \phi(D_k+1) \vert W_x \right] - E \left[ \phi(V_x+1) \vert W_x \right] \big\vert  \,  \tilde{p}_s(W_x) \right]\, dx\nn\\
& \leq &  s ||\Delta \phi||_{\infty} \int_S E^{W_x}\left[ E\left[ \vert  D_k - V_x \vert   \big\vert W_x \right]  \,  \tilde{p}_s(W_x) \right]\, dx, 
\lab{eq:UpperBound1}
\end{eqnarray}
where in the last step we have used the fact that $|\phi(i)-\phi(j)| \leq ||\Delta \phi||_{\infty} | i - j |$ for all $i, j \in \bN\cup \{0\}$. Given $A\subset \bN \cup \{0\}$ choose $\phi : \bN \cup \{0\} \to \bR$ such that  for each $i \in \bN \cup \{0\}$
\begin{equation}
1_A(i) - 1_A(Z_{\nu}) = \nu \phi(i+1) - i \phi(i)  \quad \mbox{and} \quad \phi(0)=0.
\lab{eq:g_to_indicator}
\end{equation}
By Lemma 1.1.1 of \citep{BarbourHJ} $\phi$ is bounded and satisfies $||\Delta\phi||_{\infty} \leq 1 \wedge \nu^{-1}$. (\ref{eq:stein_bound}) now follows by taking expectations in (\ref{eq:g_to_indicator}) and using (\ref{eq:UpperBound1}). \qed

To use Theorem~\ref{thm:SteinsMethod_PC} we need to construct a random variables $V_x$ coupled with $D_k$ for any $x \in S$. The heart of the problem lies in estimating an uniform upper bound for $E\left[\vert D_k - V_x\vert \big\vert W_x \right]$. The difficulty in the computations are due to the presence of random weights at each vertex and the fact that the connection function can be arbitrarily close to one. In \citep{Penrose2018} the connection function is assumed to be uniformly bounded away from one and the inhomogeneity in the graph arises from the non-uniform intensity and the fact that the (non-random) connection function is location dependent.
\subsection{Proof of Theorem~\ref{thm:poisson_conv}}
Let $x$ be a point in $S$ and $W_x$ be an independent random variable with the same distribution as the distribution of the weights. 
For $k \geq 1$ consider an extra sequence of independent and identically distributed points $Y_1, Y_2, \ldots , Y_k \in S$ with associated independent weights $W_1, \ldots ,W_k$ such that 
\[
P(Y_i \in dy |W_i) = \f{g_s\left(x,W_x; y, W_i\right) dy}{\int_{S} E^{W_i}\left[g_s\left(x,W_x; w, W_i\right)\right] dw}.
\]
Denote the graph $G(\cP_s \cup \{x\} \cup \{Y_1, Y_2, \ldots , Y_k\}, r_s)$ by $G_s$.
Let $\cP_{s, x}:=\{X \in \cP_s : X \mbox{ is neighbour of } x \mbox{ in } G_s \}$  and $\cP_s^x:=\cP_s \setminus \cP_{s, x}$. Given the weight $W_x$ at $x$ the two Poisson point processes $\cP_{s, x}$ and $\cP_s^x$  are independent. Note that 
%Given $W_x$ we define the random variables $U_x$ and $V_x$ as follows.
%
\[ D_k = D_k(\tilde{\cP}_s ) :=\#\{X \in \cP_s : \deg(X)=k \mbox{ in the graph } G_s \mbox{ induced by } \cP_s\}. \]
%

%Observe that $U_x$ does not depend on $W_x$ and has the same distribution as $D_k$. 
%
We construct a coupled point process $\cP^* \subset \cP_s \cup \{ Y_1,Y_2,\cdots, Y_k \}$. If $|\cP_{s, x}|>k$, then discard $|\cP_{s, x}|-k$ many points from $\cP_{s, x}$ chosen uniformly at random and call the resulting collection of points as $\cP^*$. If $k>0$ and $|\cP_{s, x}| < k $, then take $\cP^*=\cP_s \cup \{Y_1,Y_2,\cdots, Y_{ k -|\cP_{s, x}|}\}$.  Let $\cY_x :=\{Y_1,Y_2,\cdots, Y_{ k -|\cP_{s, x}|}\}$.
Denote the subgraph of $G_s$ induced by $\cP^* \cup \{x\}$ by $G_s^*$. Observe that $\deg(x)=k$ in $G_s^*$.  Let 
\[V_x   :=\#\{X \in \cP^* : \deg(X)=k \mbox{ in the graph } G_s^* \}. \]
Given $W_x$, the random variable $V_x$ has the same distribution as $D_k(\tilde{\cP}_s \cup \{(x, W_x)\}) - 1$  conditioned on the event $\{f(x, \cP_s)=1\}$.

If $|\cP_{s, x}|>k$ we can write $|D_k-V_x| \leq U'_x+V'_x$, where $U'_x$ is the number of vertices $y\in \cP_{s, x}$ such that $y$ is connected to $k$ points in $\cP_s$ and $V'_x$ is the number of pairs of vertices $(y, z)$, $y \in \cP_{s, x}$, $z \in \cP_s$ with $y \neq z$, $z$ is connected to $y$ and $z$ having at most $k$ neighbours in $\cP_s^x$. 

If $|\cP_{s, x}| < k $ then $|D_k-V_x| \leq U''_x+V''_x$, where $U''_x$ is the number of vertices $y\in \cY_x$ such that $y$ is connected to at most $k$ points in $\cP_s$ and $V''_x$ is the number of pairs of vertices $(y, z)$, $y \in \cY_x$, $z \in \cP_s$ with $y \neq z$, $z$ is connected to $y$ and $z$ having atmost $k$ neighbours in $\cP_s$.

By the Theorem~\ref{thm:SteinsMethod_PC} 
\begin{eqnarray}
 (1 \wedge \nu^{-1})^{-1} d_{TV}\left( F_{D_k}, F_{Z_{\nu}} \right)  & \leq &  
  s \int_S E^{W_x}\left[ E\left[ | D_k - V_x |  \big\vert W_x \right]  \,  \tilde{p}_s(W_x) \right]   dx  \nn\\
 & = &   s \int_S E^{W_x}\left[ E\left[ | D_k - V_x |   ; \{|\cP_{s, x}|>k\} \big\vert  W_x \right]  \, \tilde{p}_s(W_x) \right]   dx  \nn\\
 & & + \; s \int_S E^{W_x}\left[ E\left[ | D_k - V_x |  ; \{|\cP_{s, x}| < k\} \big\vert  W_x \right]  \,  \tilde{p}_s(W_x) \right]  dx \nn \\
& \leq  & s  \int_S E^{W_x} \left[ E\left[U'_x+V'_x \big \vert W_x \right] \, \tilde{p}_s (W_x) \right]  dx \nn\\
& & + \; s  \int_S E^{W_x} \left[ E\left[U''_x+V''_x \big \vert W_x \right] \, \tilde{p}_s(W_x) \right]  dx,
\lab{eq:TV_I}
\end{eqnarray}
where
\begin{eqnarray}
\tilde{p}_s(w) & = &   E\left[ 1{\left\{deg(O)=k \mbox{ in } G(\tilde{\cP}_s\cup \{(O, W_o)\}, r_s)\right\}} \bigg\vert  W_o=w \right]\nn\\
& = & \f{1}{k!}\left( s \int_{ S} E^{W_y}\left[g_s\left(O,w; y, W_y\right)\right] dy \right)^k \exp{\left(- s \int_{ S} E^{W_y}\left[g_s\left(O,w; y, W_y\right)\right] dy\right)}. 
\lab{eq:p}
\end{eqnarray}
By the standard change of variables 
\begin{eqnarray}
\tilde{p}_s(w) & = & \f{1}{k!}\left( s r_s^d \int_{ r_s^{-1} S} E^{W_y}\left[\tilde{g}\left(O,w; y, W_y\right)\right] dy \right)^k e^{- s r_s^d \int_{ r_s^{-1} S} E^{W_y}\left[\tilde{g}\left(O,w; y, W_y\right)\right] dy}.
\lab{eq:p_k}
\end{eqnarray}
Using the fact that the metric $d_1$ (defined in (\ref{eqn:d1})) on $r_s^{-1}S$ is toroidal and writing 
$\tilde{g}\left(w; y, W_y\right)$ for $\tilde{g}\left(O,w; y, W_y\right)$ we obtain 
\begin{eqnarray}
\tilde{p}_s(w)
& = & \f{1}{k!} \left( s r_s^d\int_{r_s^{-1} S} E^{W_y}\left[ \tilde{g}\left(w; y, W_y\right)\right] dy \right)^k \exp{\left(- s r_s^d\int_{r_s^{-1} S} E^{W_y}\left[ \tilde{g}\left(w; y, W_y\right)\right] dy \right)} \nn \\
& \leq & \f{1}{k!} \left(c_0 s r_s^d w^{\f{d}{\al}} \right)^k \exp{\left(- s r_s^d\Lam_s(w) \right)},
\lab{eq:ps(x)}
\end{eqnarray}
where the last inequality in (\ref{eq:ps(x)}) follows from the bounds in (\ref{eq:in_out_circle}).
By Theorem~\ref{thm:Ex_poisson_conv} we have $E[D_k] \to e^{-\xi}$ as $s \to \infty$. Thus it suffices to show that the right hand side of (\ref{eq:TV_I}) converges to zero. To this end we compute $E \left[ U'_x \big \vert W_x \right]$,  $E \left[ V'_x \big \vert W_x \right] $, $E \left[ U''_x \big \vert W_x \right]$ and $E \left[ V''_x \big \vert W_x \right]$. For the case $k=0$, $U''_x \equiv 0, V''_x \equiv 0$ and thus the conditions required for the theorem to hold will be determined only by the first term on the right in (\ref{eq:TV_I}).

Applying the Campbell-Mecke formula we obtain
\begin{eqnarray}
E[U'_x \big \vert W_x ]    & = &  s  \int_{S}   E^{W_y}\left[ g_s\left(x,W_x; y, W_y\right) \Phi_s (y;W_y, k) \right] dy, 
\lab{eq:EU'_x_1A}
\end{eqnarray}
where
\begin{equation}
\Phi_s (y; W_y, k):= \f{1}{k !} \left(  s \int_{S}  E^{W_z}\left[g_s\left(y,W_y; z, W_z\right) \right] dz \right)^k e^{-  s \int_{S}  E^{W_z}\left[ g_s\left(y,W_y; z, W_z\right) \right] dz}.
\lab{eq:Phi_s}
\end{equation}
Making change of variables $r_s^{-1}y \to y$ and $r_s^{-1}z \to z$ and using the fact that the metric is toroidal we obtain
\begin{eqnarray}
E[U'_x \big \vert W_x ]    & = &  s r_s^d \int_{r_s^{-1}S}   E^{W_y}\left[ \tilde{g}\left(W_x; y, W_y\right)\bar{\Phi}_s (W_y, k) \right] dy, 
\lab{eq:EU'_x_1B}
\end{eqnarray}
where
\begin{equation}
\bar{\Phi}_s (W_y, k):= \f{1}{k !} \left(  s r_s^d\int_{r_s^{-1}S}  E^{W_z}\left[\tilde{g}\left(W_y; z, W_z\right) \right] dz \right)^k e^{-  s r_s^d \int_{r_s^{-1}S}  E^{W_z}\left[ \tilde{g}\left(W_y; z, W_z\right) \right] dz}.
\lab{eq:EU'_x_1C}
\end{equation}
Bounding the expression in (\ref{eq:EU'_x_1C}) using (\ref{eq:in_out_circle}) and substituting in (\ref{eq:EU'_x_1B}) yields 
\begin{eqnarray}
E\left[U'_x \big \vert W_x \right]
& \leq  &  \f{s r_s^d}{k!}   \int_{\bR^d}   E^{W_y}\left[ \left( 1 - \exp\left( - \f{\eta W_xW_y}{|y|^\al} \right) \right) \left(c_0 s r_s^d W_y^{\f{d}{\al}}\right)^k e^{-  s r_s^d \Lam_s(W_y)  }\right] dy\nn\\
& = &  \f{s r_s^d}{k!} \int_{\bR^d}  \int_{1}^{\infty}  \left(1-\exp\left(-\f{\eta W_x w}{|y|^\al}\right)\right) \left(c_0 s r_s^d w^{\f{d}{\al}}\right)^k e^{-  s r_s^d \Lam_s(w)}  \be w^{-\be-1} dw\,dy.
\lab{eq:EU'_x_3Aa}
\end{eqnarray}
Interchanging the integrals by Fubini's theorem we have
\begin{eqnarray}
E\left[U'_x \big \vert W_x \right]
& \leq &  \f{s r_s^d}{k!}  \int_{1}^{\infty} \int_{\bR^d}  \left(1-\exp\left(-\f{\eta W_x w}{|y|^\al}\right)\right) \left(c_0 s r_s^d w^{\f{d}{\al}}\right)^k e^{-  s r_s^d \Lam_s(w)}  \be w^{-\be-1} dy\,dw.
\lab{eq:EU'_x_3A}
\end{eqnarray}
Switching to polar cordinates in the inner integral in (\ref{eq:EU'_x_3A}) and then making a change of variable $t=r^{-\al}$ we obtain
\begin{eqnarray}
\int_{\bR^d}    \left(1-\exp\left(-\f{\eta W_x w}{|y|^\al}\right)\right) dy
& = &    \int_{0}^{2\pi}  \int_{0}^{\infty}  r^{d-1}\left(1-\exp\left(-\f{\eta W_x w}{r^\al}\right)\right) dr \, d\th\nn\\
& = &    2\pi  \int_{0}^{\infty}  r^{d-1}\left(1-\exp\left(-\f{\eta W_x w}{r^\al}\right)\right) dr\nn\\
& = &     \f{2\pi}{\al}  \int_{0}^{\infty}  t^{-\f{d}{\al}-1}\left(1-e^{-\eta W_x w t}\right) dt.
\lab{eq:EU'_x_3A1}
\end{eqnarray}
Using integration by parts the integral in the last expression in (\ref{eq:EU'_x_3A1}) equals
\begin{equation}
%\int_{0}^{\infty}t^{-\f{d}{\al}-1}\Big(1- e^{-\eta W_x w t} \Big) \,dt 
%
 -\f{\al}{d} \left(1-e^{-\eta W_x w t} \right) t^{-\f{d}{\al}} \bigg\vert_{0}^{\infty} +   \f{\al}{d} (\eta W_x w) \int_{0}^{\infty} t^{-\f{d}{\al}}  e^{-\eta W_x w t} \,dt.
\lab{eq:EU'_x_3C}
\end{equation} 
Since $\al>d$, the first term in (\ref{eq:EU'_x_3C}) equals zero while the second term evaluates to  
\begin{equation}
%\int_{0}^{\infty}t^{-\f{d}{\al}-1}\Big(1- e^{-\eta W_x w t} \Big) \,dt 
%
\f{\al}{d}(\eta W_x w)^{\f{d}{\al}} \int_{0}^{\infty}  t^{-\f{d}{\al}}  e^{- t} \,du =  \f{\al}{d}\left( \eta W_x w \right)^{\f{d}{\al}} \Gamma\left(1-\f{d}{\al} \right).
\lab{eq:EU'_x_3D}
\end{equation}
Substituting from (\ref{eq:EU'_x_3D}) in (\ref{eq:EU'_x_3A1}) and then the resulting expression in (\ref{eq:EU'_x_3A}) yields
\begin{eqnarray}
E\left[U'_x \big \vert W_x \right]  
& \leq & c_0 s r_s^d (\al\be-d)\f{1}{k!} W_x^{\f{d}{\al}}   \int_{1}^{\infty} w^{\f{d}{\al}} \left(c_0 s r_s^d w^{\f{d}{\al}}\right)^k e^{-  s r_s^d \Lam_s(w)}   w^{-\be-1} \,dw \nn\\
& \leq &   \f{(\al\be-d)}{\be k!} W_x^{\f{d}{\al}}  E^{W_y} \left[  \left(c_0 s r_s^d W_y^{\f{d}{\al}}\right)^{k+1}  e^{- s r_s^d \Lam_s(W_y)   } \right]     
 =   C_1 f_1(s, k+1) W_x^{\f{d}{\al}},
\lab{eq:EU'_x_3n}
\end{eqnarray}
where
\begin{equation}
f_j(s, m):= E^{W}\left[\left( c_0 s r_s^d  W^{\f{d}{\al}} \right)^{j m} e^{ -  j s r_s^d \Lam_s(W) }\right], \qquad j \geq 1.
\lab{eq:fj}
\end{equation}
Note that $f_1(s, k+1)$ does not depend on $W_x$. Using (\ref{eq:EU'_x_3n}) and (\ref{eq:ps(x)}) we obtain
\begin{eqnarray}
\lefteqn{s r_s^d  \int_{r_s^{-1}S} E^{W_x}\left[ E\left[U'_x \big \vert W_x \right]  \, \tilde{p}_s(W_x) \right]\, dx}\nn\\ 
& \leq &  C_3 r_s^{-d} \, f_1(s, k+1)   E^{W_x} \left[\left(c_0 s r_s^d W_x^{\f{d}{\al}}\right)\left(c_0 s r_s^d W_x^{\f{d}{\al}}\right)^k   \exp\left(- s r_s^d \Lam_s(W_x)  \right) \right] \nn\\
& = &   C_3  r_s^{-d}  \,  f_1(s, k+1)^2.
%\leq C_3 r_s^{-d} \left(\f{c_1}{s}\right)^2,
%
\lab{eq:EU'_x_11nB}
\end{eqnarray}
Since $\al>d$, by Lemma~\ref{lem:equal_limit_U} with $j=1$ and $m=k+1$, $f_1(s, k+1) \leq \f{c_1}{s}$, for all $s$ sufficiently large. Hence
\begin{eqnarray}
s r_s^d  \int_{r_s^{-1}S} E^{W_x}\left[ E\left[U'_x \big \vert W_x \right]  \, \tilde{p}_s(W_x) \right]\, dx & \leq & C_4  r_s^{-d} s^{-2} \to 0 \mbox{ as } s\to \infty.
\lab{eq:EU'_x_14n}
\end{eqnarray}
Note that for the assertion in (\ref{eq:EU'_x_14n}) to hold, we only required that $\al > d$ and $\be > 1$.

Next we show that the term involving $E \left[ V'_x \big \vert W_x \right] $ in the first term on the right in (\ref{eq:TV_I}) also converges to $0$ as $s \to \infty$. By the Campbell-Mecke formula
\begin{equation}
E \left[ V'_x \big \vert W_x \right]  = \sum\limits_{i=0}^{k} s^2 \int_{S} \int_{S} E^{W_y W_z} \left[ g_s\left(x,W_x; y, W_y\right) g_s\left(y,W_y; z, W_z\right)  \f{1}{i !}  f_s(x, W_x; z, W_z; i)\, \right]dz \, dy,
\lab{eq:EV'_x_0}
\end{equation}
where
\begin{eqnarray} 
f_s(x, W_x; z, W_z; i) & = & \left(  s \int_{S}  \, E^{W_w}\left[ g_s\left(z,W_z; w, W_w\right) \left( 1- g_s\left(x,W_x; w, W_w\right) \right)\, \right] dw \right)^i \nn\\
& & \times \exp{ \left(-  s \int_{S}  \, E^{W_w}\left[ g_s\left(z,W_z; w, W_w\right) \left( 1- g_s\left(x,W_x; w, W_w\right) \right)\, \right] dw \right)}. \nn
\end{eqnarray}
If the connection function $g$ were to be uniformly bounded away from one, then replacing $1 - g$ by some $\ep > 0$ in $f_s$ above simplifies the computations considerably as in \citep{Penrose2018}. Making the standard change of variables in (\ref{eq:EV'_x_0}) 
%and use the toroidal metric to replace $\tilde{g}(x,W_x,w,W_w)$ by $\tilde{g}(W_x,w,W_w) = \tilde{g}(O,W_x,w,W_w)$ etc. This 
yields
\begin{equation}
E \left[ V'_x \big \vert W_x \right] =  \sum\limits_{i=0}^{k}  \f{s^2 r_s^{2d}}{i !} \int_{r_s^{-1}S} \int_{r_s^{-1}S}  E^{W_y W_z}\left[ \tilde{g}\left(x,W_x; y, W_y\right)    \tilde{g}\left(y,W_y; z, W_z\right)h_s(x,W_x; z,W_z; i) \right] dz \, dy,
\lab{eq:EV'_x_1n}
\end{equation}
where
\begin{eqnarray} 
h_s(x,W_x; z,W_z; i) & = & \left(  sr_s^d \int_{r_s^{-1}S}  \, E^{W_w}\left[ \tilde{g}\left(z,W_z; w, W_w\right) \left( 1- \tilde{g}\left(x,W_x; w, W_w\right) \right)\, \right] dw\right)^i \nn\\
& & \times \exp{ \left(-  sr_s^d \int_{r_s^{-1}S}  \, E^{W_w}\left[ \tilde{g}\left(z,W_z; w, W_w\right) \left( 1- \tilde{g}\left(x,W_x; w, W_w\right) \right)\right] dw \right) }. \nn
\end{eqnarray}
Using the fact that the weights are larger than one and the metric toroidal we obtain the bound 
\begin{eqnarray}
h_s(x, W_x; z, W_z; i) &\leq & \left(  sr_s^d \int_{r_s^{-1}S} \tilde{g}\left(W_z; w, W_w\right)  dw\right)^i  \nn\\
& & \times \exp\left( - s r_s^d \int_{r_s^{-1}S}  E^{W_w}\left[  \left(1- e^{-\f{\eta }{d_1(z, w)^\al}} \right) e^{-\f{\eta W_xW_w}{d_1(x, w)^\al}} \right] dw\right).
\lab{eq:EV'_x_5nA}
\end{eqnarray}	
Substituting the upper bound from (\ref{eq:in_out_circle}) in the first factor on the right in (\ref{eq:EV'_x_5nA}) yields
\begin{eqnarray}
h_s(x, W_x; z, W_z; i) \!\! & \leq & \!\!   \left( c_0 sr_s^d W_z^{\f{d}{\al}}\right)^i \exp\left( -  s r_s^d \int_{r_s^{-1}S}  \left(1-e^{-\f{\eta }{d_1(z, w)^\al}} \right)   E^{W_w}\left[  e^{-\f{\eta W_x W_w}{d_1(x, w)^\al}} \right] dw\right).
\lab{eq:EV'_x_5nB}
\end{eqnarray}
Since $\be > 1$ we have $c_1 = E[W_w] < \infty$. By the Jensen's inequality
\begin{eqnarray}
E^{W_w} \left[\exp\left(-\f{\eta W_xW_w}{d_1(x, w)^\al}\right) \right]  
\geq \exp\left(-\f{\eta  W_x c_1}{d_1(x, w)^\al}\right).
\lab{eq:EV'_x_6A}
\end{eqnarray}
Substituting from (\ref{eq:EV'_x_6A}) in (\ref{eq:EV'_x_5nB}) we obtain
\begin{equation} 
h_s(x, W_x; z, W_z; i)   \leq  \left( c_0 sr_s^d W_z^{\f{d}{\al}}\right)^i \exp\left( - s r_s^d\int_{r_s^{-1}S}    \left(1-e^{-\f{\eta }{d_1(z, w)^\al}} \right) e^{-\f{\eta c_1 W_x }{d_1(x, w)^\al}}  dw \right).
\lab{eq:EV'_x_7A}
\end{equation}
For $z, x \in r_s^{-1}S$ and $T\equiv T(s)$ to be specified later, define
\[ 
D(z, x, s ) := r_s^{-1}S \cap B(x, T^{\f{1}{\al}})^c \cap \left\{B(z, (2T)^{\f{1}{\al}})\setminus B(z, T^{\f{1}{\al}})\right\},
\]
where $B(x,r)$ is the ball centered at $x$ and radius $r$ with respect to the toroidal metric on $r_s^{-1}S$. Observe that for $w \in D(z,x,s)$, $d_1(z, w)^\al \leq 2T$ and $d_1(x, w)^\al \geq T$. Further $|D(z, x, s )| \geq c_2 T^{\f{d}{\al}}$ for some constant $c_2 > 0$. Using the above observations we obtain 
\begin{eqnarray}
\int_{r_s^{-1}S }  \left(1-e^{-\f{\eta}{d_1(z, w)^\al}} \right)  e^{-\f{\eta c_1 W_x}{d_1(x, w)^\al}}  dw 
& \geq & \int_{ D(z, x, s )}  \left(1-e^{-\f{\eta}{d_1(z, w)^\al}} \right)  e^{-\f{\eta c_1 W_x}{d_1(x, w)^\al}}  dw\nn\\
& \geq &   \int_{ D(z, x, s )}  \left(1-e^{-\f{\eta}{2T}} \right)  e^{-\f{\eta c_1 W_x}{T}}  dw \nn\\
&\geq &  \left(1-e^{-\f{\eta}{2T}}\right)  e^{-\f{\eta c_1 W_x}{T}}  c_2 T^{\f{d}{\al}}.
\lab{eq:EV'_x_8A}
\end{eqnarray}
Substituting from (\ref{eq:EV'_x_8A}) in (\ref{eq:EV'_x_7A}) we obtain 
\begin{eqnarray}
h_s(x, W_x; z, W_z; i)  & \leq & \left( c_0 sr_s^d W_z^{\f{d}{\al}}\right)^i  \exp\left( -  c_3 s r_s^d \left(1-\exp\left(-\f{\eta}{2T}\right) \right)  \exp\left(-\f{\eta c_1 W_x}{T}\right)   T^{\f{d}{\al}}\right) \nn\\
& =  &  \left( c_0 sr_s^d W_z^{\f{d}{\al}}\right)^i  \phi_s(W_x, T) \mbox{ (say)}.
\lab{eq:EV'_x_9n}
\end{eqnarray}
Since the bound in (\ref{eq:EV'_x_9n}) does not depend on $x,z$ we can use the toroidal metric and substitute in (\ref{eq:EV'_x_1n}) to obtain 
\begin{eqnarray}
E \left[ V'_x \big \vert W_x \right] 
& \leq &    \sum\limits_{i=0}^{k}  \f{s^2 r_s^{2d}}{i !} \int_{r_s^{-1}S} \int_{r_s^{-1}S}  E^{W_y W_z}\left[ \tilde{g}\left(W_x; y, W_y\right)    \tilde{g}\left(W_y; z, W_z\right)\left( c_0 sr_s^d W_z^{\f{d}{\al}}\right)^i  \phi_s(W_x, T)  \right] dz \, dy\nn\\
& \leq &  C_0 \left(s r_s^d\right)^{2+k} \phi_s(W_x, T)   \int_{r_s^{-1}S}   E^{W_y}\left[ \tilde{g}\left(W_x; y, W_y\right)  \int_{r_s^{-1}S} E^{W_z}\left[ \tilde{g}\left(W_y; z, W_z\right) W_z^{\f{kd}{\al}}  \right] dz \right] dy,
\lab{eq:EV'_x_1n1}
\end{eqnarray}
for sufficiently large $s$. Since $\al\be>(k+1)d$, proceeding as we did in (\ref{eq:EU'_x_3A})--(\ref{eq:EU'_x_3n}) yields
\begin{equation}
\int_{r_s^{-1}S} E^{W_z}\left[ \tilde{g}\left(W_y; z, W_z\right) W_z^{\f{kd}{\al}}  \right] dz 
\leq (\al\be- d) c_0 W_y^{\f{d}{\al}} E^{W_z}\left[W_z^{\f{(k+1)d}{\al}}\right]  
= \f{\al\be- d}{\al\be- (k+1)d} c_0 W_y^{\f{d}{\al}}.
\lab{eq:EV'_x_1n2}
\end{equation}
Substituting from (\ref{eq:EV'_x_1n2}) in (\ref{eq:EV'_x_1n1}) we have
\begin{eqnarray}
E \left[ V'_x \big \vert W_x \right] 
& \leq  &  C_1 \left(s r_s^d\right)^{2+k} \phi_s(W_x, T)   \int_{r_s^{-1}S}   E^{W_y}\left[ \tilde{g}\left(W_x; y, W_y\right) W_y^{\f{d}{\al}}  \right] dy\nn\\
& \leq  &  C_2 \left(s r_s^d\right)^{2+k} \phi_s(W_x, T)    c_0 W_x^{\f{d}{\al}} E^{W_y}\left[W_y^{\f{2d}{\al}}\right] \nn\\
& =  &  C_3 \left(s r_s^d\right)^{2+k} \phi_s(W_x, T)    W_x^{\f{d}{\al}}, 
\lab{eq:EV'_x_1n3}
\end{eqnarray}
where we use the fact that $E^{W_y}\left[W_y^{\f{2d}{\al}}\right] < \infty$ since $\al\be > 2d$. Using (\ref{eq:EV'_x_1n3}) and (\ref{eq:ps(x)})  we obtain
\begin{eqnarray}
%
%\lefteqn{ 
s r_s^d \int_{r_s^{-1}S}   E^{W_x} \left[ E \left[ V'_x \big\vert W_x   \right] \, \tilde{p}_s(W_x)\right]\, dx %\nn\\
%
%& \leq & C_3  s \left( sr_s^d \right)^{k+2}  E^{W_x} \left[   \phi_s(W_x, T) W_x^{\f{d}{\al}} \f{1}{k!} \left(  sr_s^d \int_{r_s^{-1}S}  E^{W_y}\left[\tilde{g}\left(W_x; y, W_y\right)\right] dy \right)^k e^{ -  sr_s^d \int_{r_s^{-1}S}  E^{W_y}\left[\tilde{g}\left(W_x; y, W_y\right)\right] dy } \right] \nn\\
%
& \leq & C_4  s \left( sr_s^d \right)^{k+2}  E^{W_x} \left[   \phi_s(W_x, T) W_x^{\f{d}{\al}} \left(c_0 s r_s^d W_x^{\f{d}{\al}}\right)^k   e^{- s r_s^d \Lam_s(W_x) }    \right]\nn\\
& = & C_5  s \left( sr_s^d \right)^{k+1}  E^{W_x} \left[   \phi_s(W_x, T)  \left(c_0 s r_s^d W_x^{\f{d}{\al}}\right)^{k+1}   e^{- s r_s^d \Lam_s(W_x)  }    \right].
\lab{eq:EV'_x_23n1}
\end{eqnarray}
An application of Cauchy-Schwarz inequality yields
\begin{eqnarray}
\qquad s r_s^d \int_{r_s^{-1}S}   E^{W_x} \left[ E \left[ V'_x \big\vert W_x   \right] \, \tilde{p}_s(W_x) \right]\, dx 
& \leq & C_5   s \left( s r_s^d \right)^{k+1}  \left(E^{W_x} \left[   \phi_s(W_x, T)^2 \right]\right)^{\f{1}{2}}  
\left( f_2(s,k+1)\right)^{\f{1}{2}}, 
\lab{eq:EV'_x_23n2}
\end{eqnarray}
where, $f_2(s,k+1)$ is defined in (\ref{eq:fj}).
We shall now compute an upper bound for $E^{W_x}\left[ \phi_s(W_x, T)^2\right]$. From (\ref{eq:EV'_x_9n}) we have for some constant $c_4$ that 
\begin{eqnarray}
\lefteqn{E^{W_x}\left[ \phi_s(W_x, T)^2\right] \; = \;
E^{W_x}\left[ \exp\left( -  c_4 s r_s^d \left( \left(1-\exp\left(-\f{\eta}{2T}\right) \right)  \exp\left(-\f{\eta c_1 W_x}{T}\right) \right) T^{\f{d}{\al}}\right)\right]}\nn\\
& = &  E^{W_x}\left[ \exp\left( -  c_4 s r_s^d \left(\left(1-\exp\left(-\f{\eta}{2T}\right) \right)  \exp\left(-\f{\eta c_1 W_x}{T}\right) \right)  T^{\f{d}{\al}}\right) ; W_x \leq T\right]\nn\\
&  & + \; E^{W_x}\left[ \exp\left( -  c_4 s r_s^d\left( \left(1-\exp\left(-\f{\eta}{2T}\right) \right)  \exp\left(-\f{\eta c_1 W_x}{T}\right) \right) T^{\f{d}{\al}}\right) ; W_x > T\right] \nn\\
& \leq &  E^{W_x}\left[ \exp\left( -  c_4 s r_s^d  \left(  \left(1-\exp\left(-\f{\eta}{2T}\right) \right)  \exp\left(-\eta c_1 \right) \right) T^{\f{d}{\al}}\right) \right]+ P\left(  W_x > T\right)\nn\\
& \leq &  \exp\left( -  c_5 s r_s^d   \left(\f{1}{T}- \f{c}{T^2} \right)  T^{\f{d}{\al}} \right) + P\left(  W_x > T\right) = R_T(s) \mbox{ (say).}
\lab{eq:EV'_x_13C}
\end{eqnarray}
Substituting from (\ref{eq:EV'_x_13C}) in (\ref{eq:EV'_x_23n2}) we obtain 
\begin{eqnarray}
s r_s^d  \int_{r_s^{-1}S} E^{W_x} \left[ E \left[ V'_x \big\vert W_x   \right] \, \tilde{p}_s(W_x) \right] \, dx 
& \leq & C_5   \left( sr_s^d \right)^{k+1}  \left(   R_T(s) \right)^{\f{1}{2}}  \left(  s^2 f_2(s,k+1)\right)^{\f{1}{2}}.
\lab{eq:EV'_x_25n}
\end{eqnarray}
By the definition of $R_T(s)$ and the fact that $s r_s^d \leq C \log s$ for $s \geq 2$, we have
\begin{equation}
\left( sr_s^d \right)^{k+1}   \left[   R_T(s) \right]^{\f{1}{2}} 
\leq    \left[ \left(\log s\right)^{2k+2}\exp\left( -  c_5 s r_s^d \left(\f{1}{T}- \f{c}{T^2} \right) T^{\f{d}{\al}} \right) +  \left(\log s\right)^{2k+2} P\left(  W_x > T\right) \right]^{\f{1}{2}}.
\lab{eq:EV'_x_11n}
\end{equation}
Since $\be > (2k+3) \left(1-\f{d}{\al}\right)$, choose $\ep>0$ such that  
\begin{equation}
\f{2k+3+\ep}{\be}\left(1 - \f{d}{\al} \right) < 1
\lab{eq:beta_n}
\end{equation}
and set $T= (\log s)^{\f{2k+3+\ep}{\be}}$. 
Recalling that $P(W_x>T)= T^{-\be}$, 
\begin{equation}
\left(\log s\right)^{2k+2} P\left(  W_x > T \right) = \f{1}{(\log s)^{1+\ep}}.
\lab{eq:T_n}
\end{equation}
Since $T \to \infty$ and $sr_s^d \geq \f{\log s}{2}$ we have for all $s$ sufficiently large
\begin{eqnarray}
\left(\log s\right)^{2k+2}\exp\left( -  c_5 s r_s^d \left(\f{1}{T}- \f{c}{T^2} \right) T^{\f{d}{\al}} \right) & \leq & \left(\log s\right)^{2k+2} \exp\left( -  c_6 s r_s^d T^{\f{d}{\al}-1} \right) \nn\\
& \leq & \left(\log s\right)^{2k+2} \exp\left( -  c_7 \log s \left((\log s)^{-\f{2k+3+\ep}{\be}\left(1 - \f{d}{\al}\right) }\right) \right).
\lab{eq:EV'_x_13_n}
\end{eqnarray}
Since $\al>d$, we have from Lemma~\ref{lem:equal_limit_U} with $j=2$ and $m=k+1$ that
\begin{equation}
s^2 f_2(s,k+1) \leq C \log s,
\lab{eq:f2logs}
\end{equation}
for all $s$ sufficiently large and some constant $C$. Thus it follows from (\ref{eq:beta_n})--(\ref{eq:f2logs}) that the expression on the right in (\ref{eq:EV'_x_25n}) converges to zero as $s \to \infty$ and hence
\begin{equation}
s r_s^d   \int_{r_s^{-1}S} E^{W_x} \left[ E \left[ V'_x \big\vert W_x   \right] \, \tilde{p}_s(W_x) \right]\, dx  \to 0 \mbox{ as } s \to \infty.
\lab{eq:2nd_term_n}
\end{equation}
In proving (\ref{eq:2nd_term_n} we have used all the three conditions mentioned in the statement of the theorem.

We shall now show that the second term on the right hand side of (\ref{eq:TV_I}) converges to 0 as $s\to \infty$. To this end we compute $E[U''_x \big \vert W_x ]$ and $E[V''_x \big \vert W_x ]$. Applying the Campbell-Mecke formula  we obtain
\begin{eqnarray}
E[U''_x \big \vert W_x ]    & = &  k \sum\limits_{i=0}^{k}  \int_{S}   E^{W_y}\left[ \f{g_s\left(x,W_x; y, W_y\right)}{\int_{S} E^{W_w}\left[ g_s\left(x,W_x; w, W_w\right)\right] dw} \; \Phi_s(y; W_y,i) \right] dy,
\lab{eq:EU''_x_00}
\end{eqnarray}
where $\Phi_s(y; W_y,i)$ is defined as in (\ref{eq:Phi_s}).
%
%\[
%\Phi_s(y; W_y,i) := \f{1}{i !}\left(  s\int_{S} %E^{W_z}\left[g_s\left(y,W_y; z, W_z\right)\right]\,dz \right)^i %\exp{ \left(-  s\int_{S} E^{W_z}\left[g_s\left(y,W_y; z, W_z\right) %\right]\,dz \right)}.
%\]
%
Making standard change of variables $r_s^{-1}y \to y$, $ r_s^{-1}z \to z$ and $r_s^{-1}w \to w$ and using the toroidal metric we can write (\ref{eq:EU''_x_00}) as
\begin{eqnarray}
E[U''_x \big \vert W_x ]  & = &    k \sum\limits_{i=0}^{k}   r_s^d \int_{r_s^{-1}S}  E^{W_y}\left[ \f{\tilde{g}\left(W_x; y, W_y\right)}{r_s^d \int_{r_s^{-1}S}E^{W_w}\left[\tilde{g}\left(W_x; w, W_w\right) \right]dw} \bar{\Phi}_s(W_y,i)  \right] dy,
\lab{eq:EU''_x_00i}
\end{eqnarray}
where $\bar{\Phi}_s(W_y,i)$ is defined as in (\ref{eq:EU'_x_1C}).
%
%\[
%\bar{\Phi}_s(W_y,i):=\f{1}{i !}\left(  sr_s^d \int_{r_s^{-1}S}  E^{W_z}\left[\tilde{g}\left(W_y; z, W_z\right)\right]dz \right)^i \exp{ \left(-  sr_s^d \int_{r_s^{-1}S}  E^{W_z}\left[\tilde{g}\left(W_y; z, W_z\right)\right] dz \right)}.
%\]
%	
Using the fact that the weights are larger than one in (\ref{eq:EU''_x_00i}) we have
\begin{equation}
E[U''_x \big \vert W_x ]    
\leq  k \sum\limits_{i=0}^{k}    \int_{r_s^{-1}S}  E^{W_y}\left[ \f{\tilde{g}\left(W_x; y, W_y\right)}{ \int_{r_s^{-1}S}\left[ 1 - \exp\left( - \f{\eta }{|w|^\al} \right) \right] dw} \bar{\Phi}_s(W_y,i) \right] dy.
\lab{eq:EU''_x_0na}
\end{equation}
Since $\al > d$ 
\begin{equation}
\int_{r_s^{-1}S}\left[ 1 - \exp\left( - \f{\eta }{|w|^\al} \right) \right] dw \to \int_{\bR^d}\left[ 1 - \exp\left( - \f{\eta }{|w|^\al} \right) \right] dw < \infty,
\lab{eq:EU''_x_0nc}
\end{equation}
and consequently we obtain the bound
\begin{eqnarray}
E[U''_x \big \vert W_x ] 
&\leq & C_1  \sum\limits_{i=0}^{k}   \int_{r_s^{-1}S}  E^{W_y}\left[\tilde{g}\left(W_x; y, W_y\right) \bar{\Phi}_s(W_y,i) \right] dy.
\lab{eq:EU''_x_0n}
\end{eqnarray}
Using the bounds from (\ref{eq:in_out_circle}) and the fact that the weights are larger than one we get
\begin{eqnarray}
E[U''_x \big \vert W_x ] 
&\leq &    C_2   \sum\limits_{i=0}^{k}  \int_{r_s^{-1}S}  E^{W_y}\left[\tilde{g}\left(W_x; y, W_y\right) \left(c_0 s r_s^d W_y^{\f{d}{\al}}\right)^i e^{-s r_s^d \Lam_s(W_y) } \right] dy\nn\\
&\leq &    C_3 \int_{\bR^d}  E^{W_y}\left[\tilde{g}\left(W_x; y, W_y\right) \left(c_0 s r_s^d W_y^{\f{d}{\al}}\right)^k e^{-s r_s^d \Lam_s(W_y) } \right] dy,
\lab{eq:EU''_x_0nb}
\end{eqnarray}
for all $s$ sufficiently large. Comparing (\ref{eq:EU''_x_0nb}) with the first inequality in (\ref{eq:EU'_x_3Aa}) and proceeding as we did to derive (\ref{eq:EU'_x_14n}) we obtain
\begin{eqnarray}
s r_s^d  \int_{r_s^{-1}S} E^{W_x}\left[ E\left[U''_x \big \vert W_x \right]  \tilde{p}_s(W_x) \right]\, dx 
%& \leq &   C_5 r_s^{-d} s^{-2} 
& \to & 0 \mbox{  as } s\to \infty.
\lab{eq:EU''_x_14n}
\end{eqnarray}
It remains to show that the term involving $ E \left[ V''_x \big \vert W_x \right] $ in the second term on the right in (\ref{eq:TV_I}) also converges to $0$ as $s \to \infty$. 
By the Campbell-Mecke formula and union bound 
\begin{eqnarray}
E \left[ V''_x \big \vert W_x \right] 
& \leq &  k  \sum_{i=0}^{k} \f{1}{i!} s  \int_{S}  \int_{S} E^{W_y W_z}\left[ \f{g_s\left(x,W_x; y, W_y\right)}{\int_{S}  E^{W_w}\left[g_s\left(x,W_x; w, W_w\right)\right] dw}     g_s\left(y,W_y; z, W_z\right)\Phi_s (z; W_z, i) \right] dz\, dy,\nn\\
\lab{eq:EV''_x_1}
\end{eqnarray}
where $\Phi_s(z; W_z,i)$ is defined as in (\ref{eq:Phi_s}).
%
%\[ 
%\Phi_s (z; W_z,  i):= \f{1}{i !} \left(  s \int_{S}  E^{W_w}\left[g_s\left(z,W_z; w, W_w\right) \right] dw \right)^i e^{-  s \int_{S}  E^{W_w}\left[ g_s\left(z,W_z; w, W_w\right) \right] dw}.
%\]
%
By making the standard change of variable the expression inside the sum on the right in (\ref{eq:EV''_x_1}) can be bounded by some constant times
\begin{equation}
%E \left[ V''_x \big \vert W_x \right] 
%
s r_s^{2d} \int_{r_s^{-1}S}  \int_{r_s^{-1}S} E^{W_y W_z}\left[ \f{\tilde{g}\left(W_x; y, W_y\right)}{r_s^d\int_{r_s^{-1}S}  E^{W_w}\left[\tilde{g}\left(W_x; w, W_w\right)\right] dw}    \tilde{g}\left(W_y; z, W_z\right) \bar{\Phi}_s (W_z, i) \right] dz\, dy,
\lab{eq:EV''_x_a1A}
\end{equation}
where $\bar{\Phi}_s(W_z,i)$ is defined as in (\ref{eq:EU'_x_1C}).
%
%\[
%\bar{\Phi}_s ( W_z, i) = \f{1}{i !} \left(  s r_s^d\int_{r_s^{-1}S}  E^{W_w}\left[\tilde{g}\left(W_z; w, W_w\right) \right] dw \right)^i e^{ - s r_s^d \int_{r_s^{-1}S}  E^{W_w}\left[\tilde{g}\left(W_z; w, W_w\right)\right] dw }.
%\]
%
Using the fact that the weights are larger than one the expression in (\ref{eq:EV''_x_a1A}) is bounded by
%
% E \left[ V''_x \big \vert W_x \right] 
%
% &\leq &  k  \sum_{i=0}^{k} \f{1}{i!} 
%
\[ s r_s^{d} \int_{r_s^{-1}S}  \int_{r_s^{-1}S} E^{W_y W_z} \left[ \f{\tilde{g}\left(W_x; y, W_y\right)}{\int_{r_s^{-1}S}\left[ 1 - \exp\left( - \f{\eta }{|w|^\al} \right)\right]dw} \tilde{g}\left(W_y; z, W_z\right) \bar{\Phi}_s (W_z, i)  \right] dz\, dy \]
%
% &\leq &  C_1  \sum_{i=0}^{k} \f{1}{i!} 
%
\begin{equation}
\leq \; C_0 s r_s^d  \int_{r_s^{-1}S} E^{W_y}\left[ \tilde{g}\left(W_x; y, W_y\right) \int_{r_s^{-1}S} E^{W_z}\left[  \tilde{g}\left(W_y; z, W_z\right) \bar{\Phi}_s (W_z, i) \right] dz \right] dy,
\lab{eq:EV''_x_a1B}
\end{equation}
where the last inequality follows from (\ref{eq:EU''_x_0nc}).
%by following the same arguments as in (\ref{eq:EU''_x_0na})--(\ref{eq:EU''_x_0n}), since $\al>d$.
%
Using (\ref{eq:in_out_circle}) and the fact that weights are larger than one, we bound the inner integral in (\ref{eq:EV''_x_a1B}) as follows. 
\begin{eqnarray}
\int_{r_s^{-1}S}  E^{W_z}\left[  \tilde{g}\left(W_y; z, W_z\right) \bar{\Phi}_s (W_z, i) \right] dz 
&\leq & C_1 \int_{r_s^{-1}S}  E^{W_z}\left[  \tilde{g}\left(W_y; z, W_z\right)  \left(c_0 s r_s^d W_z^{\f{d}{\al}}\right)^i e^{- s r_s^d \Lam_s(W_z) } \right] dz\nn\\
&\leq & C_1 \int_{r_s^{-1}S}  E^{W_z}\left[  \tilde{g}\left(W_y; z, W_z\right)  \left(c_0 s r_s^d W_z^{\f{d}{\al}}\right)^k e^{- s r_s^d \Lam_s(W_z) } \right] dz,
\lab{eq:Psi_2n}
\end{eqnarray}
for sufficiently large $s$. 
Comparing (\ref{eq:Psi_2n}) with the first inequality in (\ref{eq:EU'_x_3Aa}) and proceeding as we did to derive (\ref{eq:EU'_x_3n}) we obtain 
\begin{eqnarray}
\int_{r_s^{-1}S}  E^{W_z}\left[  \tilde{g}\left(W_y; z, W_z\right)  \left(c_0 s r_s^d W_z^{\f{d}{\al}}\right)^k e^{- s r_s^d \Lam_s(W_z) } \right] dz
&\leq & C_2 f_1(s,k+1) c_0 W_y^{\f{d}{\al}},
\lab{eq:Psi_2nA}
\end{eqnarray}
where  $f_1(s,k+1)$ is defined in (\ref{eq:fj}). Substituting from (\ref{eq:Psi_2n}) and (\ref{eq:Psi_2nA}) in (\ref{eq:EV''_x_a1B}) we obtain 
\begin{eqnarray}
E \left[ V''_x \big \vert W_x \right] 
&\leq &  C_3  s r_s^d \int_{r_s^{-1}S}   E^{W_y}\left[\tilde{g}\left(W_x; y, W_y\right) f_1(s,k+1) c_0 W_y^{\f{d}{\al}}  \right] dy\nn\\
& = &  C_3 f_1(s,k+1) c_0 s r_s^d \int_{r_s^{-1}S}   E^{W_y}\left[ \tilde{g}\left(W_x; y, W_y\right) W_y^{\f{d}{\al}}  \right] dy\nn\\
& \leq & C_4  s r_s^d f_1(s,k+1) W_x^{\f{d}{\al}} E^{W_y}\left[W_y^{\f{2d}{\al}}\right]= C_5  s r_s^d f_1(s,k+1) W_x^{\f{d}{\al}},
\lab{eq:EV''_x_an1}
\end{eqnarray}
since $\al\be>2d$. From (\ref{eq:EV''_x_an1}) and (\ref{eq:ps(x)}) we get
\begin{eqnarray}
%
%\lefteqn{
s r_s^d \! \int_{r_s^{-1}S} E^{W_x}\left[ E\left[V''_x \big \vert W_x \right] \tilde{p}_s(W_x) \right] dx %}\nn\\
& \leq  &   C_5 s^2  r_s^{d} f_1(s, k+1)  E^{W_x}\left[  W_x^{\f{d}{\al}} \left(c_0 s r_s^d W_x^{\f{d}{\al}}\right)^{k}  e^{- s r_s^d \Lam_s(W_x)} \right]\nn\\
& =  &   C_6 s   f_1(s, k+1)  E^{W_x}\left[ c_0 s r_s^{d} W_x^{\f{d}{\al}} \left(c_0 s r_s^d W_x^{\f{d}{\al}}\right)^{k}  e^{- s r_s^d \Lam_s(W_x)} \right]\nn\\
& = & C_6 s f_1(s, k+1)^2.
\lab{eq:EV''_x_4n}
\end{eqnarray}
Again by Lemma~\ref{lem:equal_limit_U} with $j=1$ and $m=k+1$ we have $f_1(s, k+1 ) \leq \f{c_1}{s}$, for all $s$ sufficiently large. Hence
\begin{eqnarray}
s r_s^d  \int_{r_s^{-1}S} E^{W_x}\left[ E\left[V''_x \big \vert W_x \right]  \tilde{p}_s(W_x) \right]\, dx &\leq &    C_7 s^{-1} \to 0 \mbox{ as } s\to \infty. 
\lab{eq:EV''_x_9n}
\end{eqnarray}
In order to prove (\ref{eq:EV''_x_9n}) we have used the conditions that $\al>d, \be>1$ and $\al\be>2d$.
Theorem~\ref{thm:poisson_conv} now follows from (\ref{eq:EU'_x_14n}), (\ref{eq:2nd_term_n}), (\ref{eq:EU''_x_14n}) and (\ref{eq:EV''_x_9n}). \qed
\subsection{Proof of Theorem~\ref{thm:connectivity}} 

Let $\hat{r}_s(\gamma), \kappa, T, Q$ be as defined in (\ref{eq:def_kappa}), (\ref{eq:def_TQ}) and $\rho$ satisfy $Q(\rho) = 1$. Fix $\gamma > \rho$. Since $Q$ is continuous and increasing (see remark below (\ref{eq:def_TQ})), we can and do choose $b>1$ such that $\gamma T(\f{\gamma}{b})>1$. Let $\tilde{r}_s(b)^d :=\f{b \log s}{\th_d s}$. 
In any graph $G$ with vertex set $V$ and edge set $E$, a one-hop path between distinct vertices $x,y \in V$ is a two-path comprising of edges $\{x,z\}, \{z,y\} \in E$  for some $z \in V$.  Let $E_s \equiv E_s(\gamma, b)$ be the event that there is a vertex $X\in \cP_s$ such that $X$ does not have a one-hop path to some vertex in $\cP_s\cap B(X, \tilde{r}_s(b))$ in the graph $G(\cP_s, \hat{r}_s(\gamma))$. The result is not altered if we restrict attention only to those vertices in $\cP_s\cap B(X, \tilde{r}_s(b))$ that do not have a direct connection to $X$.
 
If $P(E_s)\rar 0$ as $s\rar\infty$, then every vertex in the graph $G(\cP_s, \hat{r}_s(\gamma))$ is connected to a point in $\cP_s\cap B(X, \tilde{r}_s(b))$ through an one-hop path whp. The result then follows since the existence of a path to all the $ \tilde{r}_s(b)$-neighbours for some $b>1$ will imply that the graph $G(\cP_s, \hat{r}_s(\gamma))$ is connected whp ( Theorem 13.7, \citep{Penrose2003}). Thus it suffices to show that $P(E_s)\rar 0$ as $s\rar\infty$.

Let $h(x,\cP_s)$ equal one if $x$ is not connected to at least one vertex in $\cP_s \cap B(x, \tilde{r}_s(b))$ via a one hop path in the graph $G(\cP_s \cup \{x\}, \hat{r}_s(\gamma))$ and zero otherwise. By the Campbell-Mecke formula we have
\begin{equation}
P(E_s)  \leq  E\Big[\sum_{x\in\cP_s}h(x,\cP_s) \Big] = s E^o[h(O,\cP_s)].
\lab{eq:E_s}
\end{equation}
Let $H:\bR^+\rar \bR^+$ defined by $H(x)=1-x+x\log x$.  Choose $a$ large enough such that $bH\big(\f{a}{b}\big)>1$. Define the event $F_s :=\left\{ \cP_s \left( B(O,\tilde{r}_s(b)) \right) < a\log s \right\}$. From (\ref{eq:E_s}) we have
\begin{equation}
P(E_s)  \leq s E^o[h(O,\cP_s); F_s] + s P(F_s^c).
\lab{eq:E_s1}
\end{equation}
By the Chernoff bound (see Lemma 1.2 of \citep{Penrose2003}), 
\begin{equation}
s  P(F_s^c) =  s  P\left( \cP_s \left( B(O,\tilde{r}_s(b)) \right) \geq a\log s \right)  \leq   s \exp\left(-s\th \tilde{r}_s(b)^d H \left( \f{a\log s}{s\th \tilde{r}_s(b)^d} \right) \right) 
=  \f{1}{s^{bH\left( \f{a}{b} \right)-1}} \to 0, 
\lab{eq:F_s}
\end{equation}
as  $s \to \infty$, since $bH \left(\f{a}{b} \right) > 1$. By conditioning on the number of points of $\cP_s$ in the ball $B(O, \tilde{r}_s(b))$ the first term on the right in~(\ref{eq:E_s1}) can be bounded by
\begin{equation}
E^o[h(O,\cP_s); F_s] \leq  a\log s \; P^o(A_s),
\lab{eq:A_s_H_s}
\end{equation}
where $A_s$ is the event that $O$ is not connected to a point chosen uniformly at random in $ B(O, \tilde{r}_s(b))$ via a one-hop path in the graph  $G(\cP_s, \hat{r}_s(\gamma))$. By the thinning theorem (see Proposition 5.5, Theorem 5.8, \citep{LastPenrose})
\begin{equation}
P^o(A_s)   =  \f{1}{\th_d\tilde{r}_s(b)^d} \int_{B(O, \tilde{r}_s(b))} E^{W_0 W_y}\left[\exp\left(-s\int_{S}E^{W_z}\left[ \hat{g}_s(O,W_o; z,W_z)\, \hat{g}_s(y,W_y; z,W_z) \right] dz \right) \right] dy.
\lab{eq:A_s}
\end{equation}
Since the weights are all greater than one, the integral inside the exponential in (\ref{eq:A_s}) can be bounded from below as follows. Since $\hat{r}_s(\gamma) \to 0$ as $s \to \infty$, we have for all $s$ sufficiently large,
\begin{eqnarray}
\lefteqn{ \int_{S} E^{W_z}\left[ \left(1-\exp\left(-\f{\eta W_0W_z}{\left(\f{d(O,z)}{\hat{r}_s(\gamma)}\right)^\al}\right)\right) \left(1-\exp\left(-\f{\eta W_yW_z}{\left(\f{d(y, z)}{\hat{r}_s(\gamma)}\right)^\al}\right)\right)\right]\,dz  }\nn\\
& \geq & \int_{B(y, \hat{r}_s(\gamma))} \left(1-\exp\left(-\f{\eta }{\left(\f{d(O,z)}{\hat{r}_s(\gamma)}\right)^\al}\right)\right) E^{W_z}\left[1-\exp\left(-\f{\eta W_z}{\left(\f{d(y, z)}{\hat{r}_s(\gamma)}\right)^\al}\right)\right]\,dz.
\lab{eq:A_s_1a}
\end{eqnarray}
Using the triangle inequality in the second factor inside the integral and a change of variable the expression on the right in (\ref{eq:A_s_1a}) can be bounded from below by
\begin{equation}
%& \geq &  
\left(1-\exp\left(-\eta \bigg\vert\f{\hat{r}_s(\gamma)+\tilde{r}_s(b)}{\hat{r}_s(\gamma)}\bigg\vert^{-\al}\right)\right) \hat{r}_s(\gamma)^d \int_{B(O , 1)}  E^{W_z}\left[1-\exp\left(- \f{\eta W_z}{|z|^{\al}}\right)\right]\,dz  =  \kappa T\left(\f{\gamma}{b}\right) \hat{r}_s(\gamma)^d,
%\nn\\
%
  %& = & T\left(\f{\gamma}{b}\right) \hat{r}_s(\gamma)^d \int_{B(y, 1)} E^{W_z}\left[1-\exp\left(-\eta W_z\right)\right]\,dz \nn\\
%
% & = & \kappa T\left(\f{\gamma}{b}\right) \hat{r}_s(\gamma)^d,
%  
\lab{eq:A_s_1}
\end{equation}
since $\f{\tilde{r}_s(b)}{\hat{r}_s(\gamma)}=  \left( \f{b \kappa }{\gamma \th_d} \right)^{\f{1}{d}}$.
From (\ref{eq:A_s_H_s})--(\ref{eq:A_s_1}) we obtain 
\begin{eqnarray}
s \, E^o[h(O,\cP_s); F_s]
& \leq & a  s \log s \exp\left( - \kappa  T\left(\f{\gamma}{b}\right) s \hat{r}_s(\gamma)^d \right) \nn\\
& = & a  s^{- \gamma T\left(\f{\gamma}{b}\right)+1}  \log s \to 0,  
\lab{eq:A_s_2}
\end{eqnarray}
as $s \to \infty$ since $ \gamma T\left(\f{\gamma}{b}\right)>1$. It follows from (\ref{eq:E_s1}), (\ref{eq:F_s}) and (\ref{eq:A_s_2}) that $P(E_s)\to 0$ as $s \to \infty$. This completes the proof of Theorem~\ref{thm:connectivity} .\qed

%\bibliographystyle{plain}
%\bibliography{poisson_conv}
\end{document}